# The Generalized Marshall-Olkin-Kumaraswamy-G family of distributions

Laba Handique and Subrata Chakraborty*

Department of Statistics, Dibrugarh University

Dibrugarh-786004, Assam, India

*Corresponding Author. Email: subrata_stats@dibru.ac.in

(August 21, 2016)

**Abstract**

A new family of distribution is proposed by using Kumaraswamy-*G* (Cordeiro and de Castro, 2011) distribution as the base line distribution in the Generalized Marshal-Olkin (Jayakumar and Mathew, 2008) Construction. A number of special cases are presented. By expanding the probability density function and the survival function as infinite series the proposed family is seen as infinite mixtures of the Kumaraswamy-*G* (Cordeiro and de Castro, 2011) distribution. Density function and its series expansions for order statistics are also obtained. Order statistics, moments, moment generating function, Rényi entropy, quantile function, random sample generation, asymptotes, shapes and stochastic orderings are also investigated. The methods of parameter estimation by method of maximum likelihood and method of moment are presented. Large sample standard error and confidence intervals for the mles are also discussed. One real life application of comparative data fitting with some of the important sub models of the family and some other models is considered.

---

*Preprint*

## 1. Introduction

Generating new distributions starting with a base line distribution by adding one or more additional parameters through various mechanisms is an area of research in the filed of the probability distribution which have seen lot of work of late. The basic motivation of these works is to bring in more flexibility in the modelling different type of data generated from real life situation.

Recently, there is renewed activity in this area to propose and investigate new families of distributions. A recent review paper by Tahir *et al.* (2015) provides a detail account of various popular techniques of generating new families of univariate continuous distributions through introduction of additional parameters. Very recent contributions in this line include the Kumaraswamy Marshal-Olkin family proposed by Alizadeh *et al.*, (2015) and Marshal - Olkin Kumaraswamy-*G* family introduced by Handique and Chakraborty (2015) among others.

In this article we propose another family of continuous probability distribution by integrating the Kumaraswamy-*G* family (Cordeiro and de Castro, 2011) as the base line distribution in the Generalized Marshall Olkin Extended family (Jayakumar and Mathew, 2008). This new family referred to as the Generalized Marshall-Olkin Kumaraswamy-*G* ($GMOKw-G$) family of distribution is investigated for some its general properties, parameter estimation and real life applications.

The rest of this article is organized in seven sections. In section 2 we briefly introduce some important characteristics of probability distributions, the generalized Marshall-Olkin (Jayakumar and Mathew, 2008) and the Kumaraswamy-*G* (Cordeiro and de Castro, 2011) family of distributions. The proposed new family is defined along with its physical basis in section 3. Important special cases of the family along with their shape and main reliability characteristics are presented in the next section. In section 5 we discuss some general results of the proposed family, while different methods of estimation of parameters along comparative data modelling example are presented in section 6. The article ends with a conclusion in section 7 followed by an appendix to derive asymptotic confidence bounds.

## 2. Some formulas and notations

Here first we list some formulas to be used in the subsequent sections of this article.

If $T$ is a continuous random variable with pdf, $f(t)$ and cdf $F(t)=P[T \le t]$, then its

Survival function (sf): $\bar{F}(t)=P[T>t]=1-F(t)$,

Hazard rate function (hrf): $h(t)=f(t)/\bar{F}(t)$,

Reverse hazard rate function (rhrf): $r(t)=f(t)/F(t)$,

Cumulative hazard rate function (chrf): $H(t)=-\log[\bar{F}(t)]$,

$(p,q,r)^{\text{th}}$ Power Weighted Moment (PWM): $\Gamma_{p,q,r} = \int_{-\infty}^{\infty} t^p [F(t)]^q [1-F(t)]^r f(t)\,dt$,

Rényi entropy: $I_R(\delta) = (1-\delta)^{-1} \log\left( \int_{-\infty}^{\infty} f(t)^{\delta}\, dt \right)$.

**2.1 Generalized Marshall-Olkin Extended (** GMOE **) family of distributions**

Jayakumar and Mathew (2008) proposed a generalization of the Marshall and Olkin (1997) family of distributions by using the Lehman second alternative (Lehmann 1953) to obtain the sf $\bar{F}^{GMO}(t)$ of the GMOE family of distributions by exponentiation the sf of MOE family of distributions as

$$\bar{F}^{GMO}(t) = \left[ \frac{\alpha \bar{F}(t)}{1-\bar{\alpha}\bar{F}(t)} \right]^{\theta}, -\infty < t < \infty; 0 < \alpha < \infty; 0 < \theta < \infty \qquad (1)$$

where $-\infty < t < \infty$, $\alpha > 0\,(\bar{\alpha} = 1-\alpha)$ and $\theta > 0$ is an additional shape parameter. When $\theta = 1$, $\bar{F}^{GMO}(t) = \bar{F}^{MO}(t)$ and for $\alpha = \theta = 1$, $\bar{F}^{GMO}(t) = \bar{F}(t)$. The cdf and pdf of the GMOE distribution are respectively

$$F^{GMO}(t) = 1 - \left[ \frac{\alpha \bar{F}(t)}{1-\bar{\alpha}\bar{F}(t)} \right]^{\theta} \qquad (2)$$

and $$f^{GMO}(t) = \theta \left[ \frac{\alpha \bar{F}(t)}{1-\bar{\alpha}\bar{F}(t)} \right]^{\theta-1} \left\{ \frac{\alpha f(t)}{[1-\bar{\alpha}\bar{F}(t)]^2} \right\} = \frac{\theta \alpha^{\theta} f(t) \bar{F}(t)^{\theta-1}}{[1-\bar{\alpha}\bar{F}(t)]^{\theta+1}} \qquad (3)$$

Reliability measures like the hrf, rhrf and chrf associated with (1) are

$$h^{GMO}(t) = \frac{f^{GMO}(t)}{\bar{F}^{GMO}(t)} = \frac{\theta \alpha^{\theta} f(t) \bar{F}(t)^{\theta-1}}{[1-\bar{\alpha}\bar{F}(t)]^{\theta+1}} \Bigg/ \left[ \frac{\alpha \bar{F}(t)}{1-\bar{\alpha}\bar{F}(t)} \right]^{\theta}$$

$$= \theta \frac{f(t)}{\bar{F}(t)} \frac{1}{1-\bar{\alpha}\bar{F}(t)} = \theta \frac{h(t)}{1-\bar{\alpha}\bar{F}(t)}$$

$$r^{GMO}(t) = \frac{f^{GMO}(t)}{F^{GMO}(t)} = \frac{\theta \alpha^{\theta} f(t) \bar{F}(t)^{\theta-1}}{[1-\bar{\alpha}\bar{F}(t)]^{\theta+1}} \Bigg/ 1 - \left[ \frac{\alpha \bar{F}(t)}{1-\bar{\alpha}\bar{F}(t)} \right]^{\theta}$$

$$= \frac{\theta \alpha^{\theta} f(t) \bar{F}(t)^{\theta-1}}{[1-\bar{\alpha}\bar{F}(t)]\,[\{1-\bar{\alpha}\bar{F}(t)\}^{\theta} - \alpha^{\theta}\bar{F}(t)^{\theta}]} = \frac{\theta \alpha^{\theta} f(t) \bar{F}(t)^{\theta-1}}{[1-\bar{\alpha}\bar{F}(t)]^{\theta+1} - \alpha^{\theta}\bar{F}(t)^{\theta}[1-\bar{\alpha}\bar{F}(t)]}$$

$$H^{GMO}(t) = -\log \left[ \frac{\alpha \bar{F}(t)}{1-\bar{\alpha}\bar{F}(t)} \right]^{\theta} = -\theta \log \left[ \frac{\alpha \bar{F}(t)}{1-\bar{\alpha}\bar{F}(t)} \right]$$

Where $f(t), F(t)$, $\bar{F}(t)$ and $h(t)$ are respectively the pdf, cdf, sf and hrf of the baseline distribution.

We denote the family of distribution with pdf (1) as $\text{GMOE}(\alpha,\theta,)$ (Jayakumar and Mathew, 2008) which for $\theta=1$, reduces to $\text{MOE}(\alpha)$ (Marshall and Olkin, 1997).

**2.2 Kumaraswamy-G $(Kw-G)$ family of distributions**

For a baseline cdf $G(t)$ with pdf $g(t)$, Cordeiro and de Castro (2011) defined $Kw-G$ distribution with cdf and pdf

$$F^{KwG}(t)=1-[1-G(t)^a]^b,\quad 0<t<\infty,\, 0<a,b<\infty \tag{4}$$

and

$$f^{KwG}(t)=a\,b\,g(t)G(t)^{a-1}[1-G(t)^a]^{b-1} \tag{5}$$

Where $t>0,\ g(t)=\frac{d}{dt}G(t)$ and $a>0, b>0$ are shape parameters in addition to those in the baseline distribution. The sf, hrf, rhrf and chrf of this distribution are respectively given by

$$\overline{F}^{KwG}(t)=1-F^{KwG}(t)=[1-G(t)^a]^b$$

$$h^{KwG}(t)=f^{KwG}(t)/\overline{F}^{KwG}(t)=a\,b\,g(t)G(t)^{a-1}[1-G(t)^a]^{b-1}/[1-G(t)^a]^b$$

$$=a\,b\,g(t)G(t)^{a-1}[1-G(t)^a]^{-1}$$

$$r^{KwG}(t)=f^{KwG}(t)/\overline{F}^{KwG}(t)=a\,b\,g(t)G(t)^{a-1}[1-G(t)^a]^{b-1}/1-[1-G(t)^a]^b$$

$$=a\,b\,g(t)G(t)^{a-1}[1-G(t)^a]^{b-1}\left\{1-[1-G(t)^a]^b\right\}^{-1}$$

and $H^{KwG}(t)=-b\log[1-G(t)^a]$.

**3. Generalized Marshall-Olkin Kumaraswamy-*G* $(GMOKw-G)$ family of distributions**

We now propose a new extension of the GMO family by considering the cdf and pdf of $Kw-G$ distribution in (4) and (5) as the $f(t)$ and $F(t)$ respectively in the GMO formulation in (3) and call it $GMOKw-G$ distribution. The resulting expression for the pdf of $GMOKw-G$ is given by

$$f^{GMOKwG}(t)=\frac{\theta\alpha^{\theta}\,a\,b\,g(t)G(t)^{a-1}[1-G(t)^a]^{b-1}[1-G(t)^a]^{b(\theta-1)}}{[1-\overline{\alpha}[1-G(t)^a]^b]^{\theta+1}}$$

$$=\frac{\theta\alpha^{\theta}\,a\,b\,g(t)G(t)^{a-1}[1-G(t)^a]^{b\theta-1}}{[1-\overline{\alpha}[1-G(t)^a]^b]^{\theta+1}}, \tag{6}$$

$$0<t<\infty,\alpha>0,\theta>0,\, 0<a,b<\infty$$

The cdf, sf, hrf, rhrf and chrf of $GMOKw-G$ distribution are respectively given by

$$F^{GMOKwG}(t)=1-\left[\frac{\alpha[1-G(t)^a]^b}{1-\overline{\alpha}[1-G(t)^a]^b}\right]^{\theta} \text{ and } \overline{F}^{GMOKwG}(t)=\left[\frac{\alpha[1-G(t)^a]^b}{1-\overline{\alpha}[1-G(t)^a]^b}\right]^{\theta} \tag{7}$$

hrf:

$$h^{GMOKwG}(t)=\frac{\theta\,a\,b\,g(t)G(t)^{a-1}[1-G(t)^a]^{-1}}{1-\overline{\alpha}[1-G(t)^a]^b} \tag{8}$$

rhrf : $r^{GMOKwG}(t) = \dfrac{\theta\alpha^{\theta}\, a b\, g(t) G(t)^{a-1}[1-G(t)^{a}]^{b\theta-1}}{[1-\overline{\alpha}[1-G(t)^{a}]^{b}]^{\theta+1} - \alpha^{\theta}[1-G(t)^{a}]^{b\theta}[1-\overline{\alpha}[1-G(t)^{a}]^{b}]}$

chrf: $H^{GMOKwG}(t) = -\theta \log\left[\dfrac{\alpha[1-G(t)^{a}]^{b}}{1-\overline{\alpha}[1-G(t)^{a}]^{b}}\right]$

For $\theta = 1$, we get back the corresponding expressions for $MOKw-G$ distribution of Handique and Chakraborty (2015), if $a = b = 1$ then $f^{GMOKwG}(t) = f^{GMO}(t)$ (Jayakumar and Mathew, 2008) and for, $\alpha = 1$, $f^{GMOKwG}(t) = f^{GKwG}(t)$.

**3.1 Genesis of the distribution**

For $i = 1, 2, ..., \theta$ where $\theta > 0$ is an integer, if $\{T_{i1}, T_{i2} ..., T_{iN}\}$ be a sequence of *i.i.d.* random variables with survival function $[1-G(t)^{a}]^{b}$, and

i. if $N$ has a geometric distribution with parameter $\alpha\,(0 < \alpha \le 1)$ independent of $T_{ij}$'s, then $\min\limits_{1\le i\le\theta}\{\min(T_{i1}, T_{i2} ..., T_{iN})\}$ is distributed as $\text{GMOKw-G}(\alpha, \theta, a, b)$ or

ii. if $N$ has a geometric distribution with parameter $1/\alpha\,(\alpha > 1)$ independent of $T_{ij}$'s, then $\min\limits_{1\le i\le\theta}\{\max(T_{i1}, T_{i2} ..., T_{iN})\}$ is distributed as $\text{GMOKw-G}(\alpha, \theta, a, b)$.

***Proof.***

Let $T_{i1}, T_{i2} ..., T_{iN}$ be a sequence of *i.i.d.* random variables with survival function $[1-G(t)^{a}]^{b}$, and suppose $N$ has a geometric distribution with parameter $p$ independent of $T_{ij}$ 's. Then it can easily shown that $W_i = \min(T_{i1}, T_{i2} ..., T_{iN})$ and $V_i = \max(T_{i1}, T_{i2} ..., T_{iN})$ are distributed as $\text{MOKw-G}(p, a, b)$ and $\text{MOKw-G}(p^{-1}, a, b)$ respectively.

i. For $0 < \alpha \le 1$, if $N$ has a geometric distribution with parameter $\alpha$ ,then

$$P[\min\{W_1, W_2, ..., W_\theta\} > t]$$

$$= P[W_1 > t]\, P[W_2 > t] ... P[W_\theta > t]$$

$$= \prod_{i=1}^{\theta} P[W_i > t] = [\overline{F}^{MOKwG}(t; \alpha, a, b)]^{\theta}$$

$$= \left[\frac{\alpha[1-G(t)^{a}]^{b}}{1-\overline{\alpha}[1-G(t)^{a}]^{b}}\right]^{\theta}$$

This is the sf of $\text{GMOKw-G}(\alpha, \theta, a, b)$.

ii. For $\alpha > 1$, if $N$ has a geometric distribution with parameter $\alpha^{-1}$, then

$$P[\min\{V_1, V_2, ..., V_\theta\} > t]$$

$$= P[V_1 > t]\, P[V_2 > t] \ldots P[V_\theta > t]$$

$$= \prod_{i=1}^{\theta} P[V_i > t] = [\overline{F}^{\,MOKwG}(t;\alpha,a,b)]^{\theta},\ [\because (\alpha^{-1})^{-1} = \alpha\,]$$

$$= \left[ \frac{\alpha [1 - G(t)^{a}]^{b}}{1 - \overline{\alpha} [1 - G(t)^{a}]^{b}} \right]^{\theta}.$$

## 2.2 Shape of the density and hazard functions

Here we have plotted the pdf and hrf of the $GMOKw-G$ for some choices of the parameters to study the variety of shapes assumed by the family.

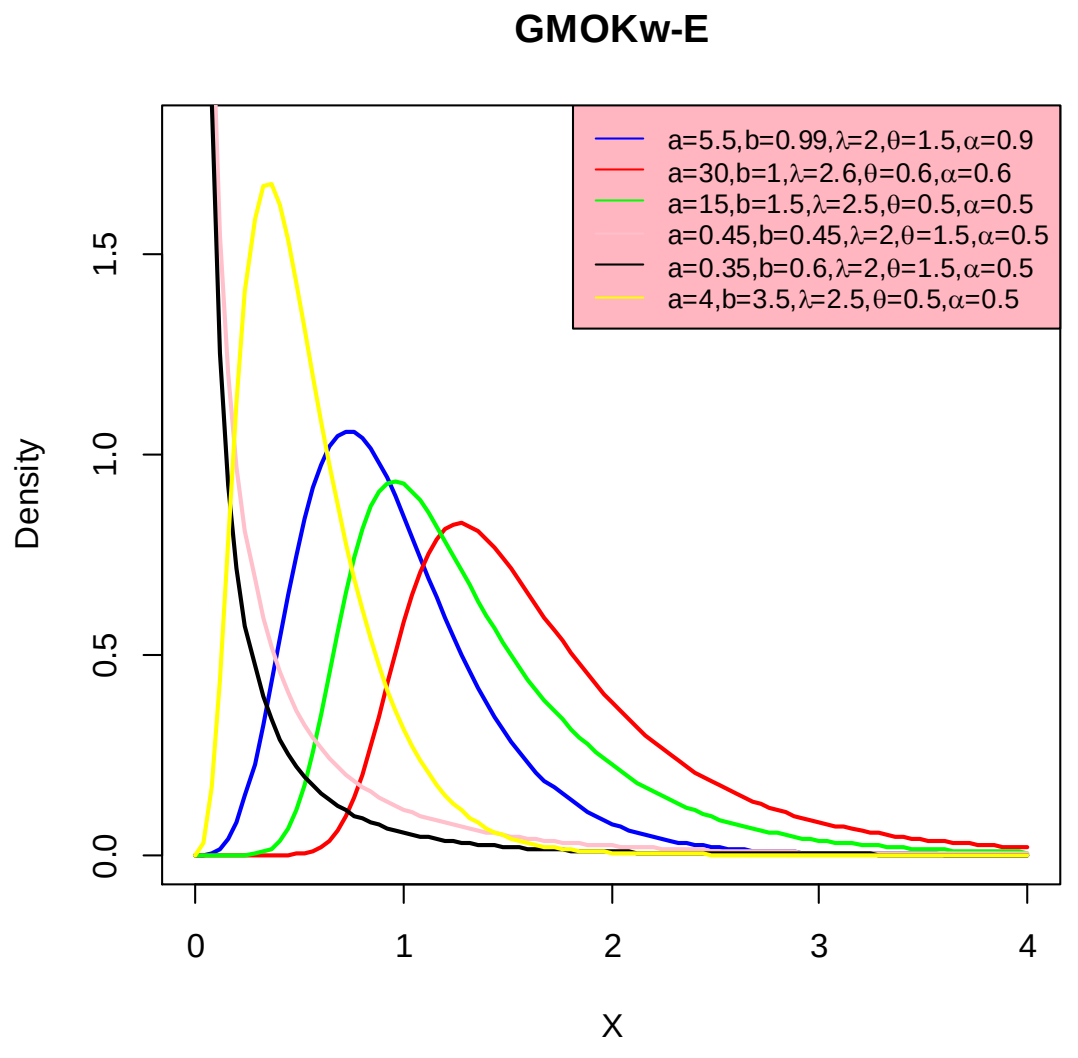

(a)

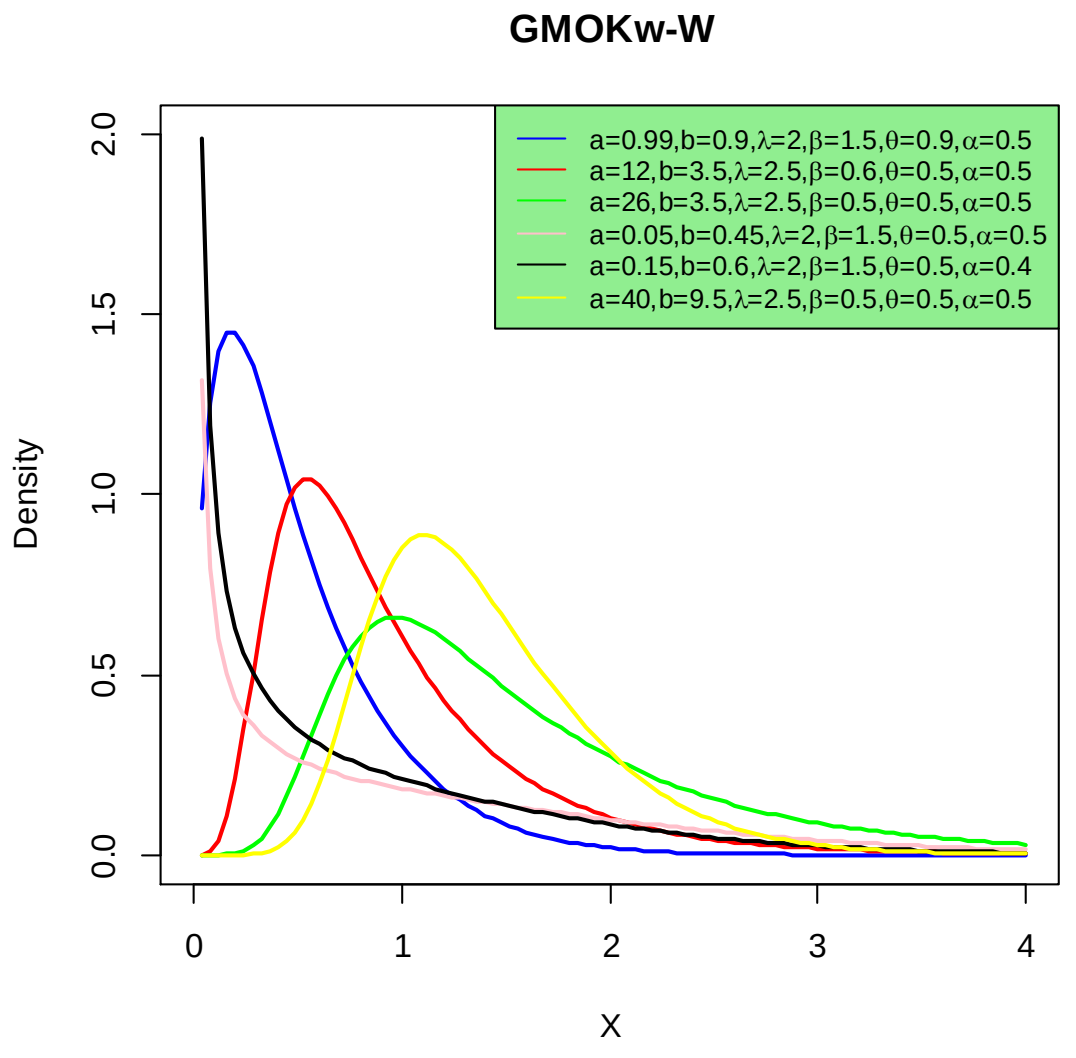

(b)

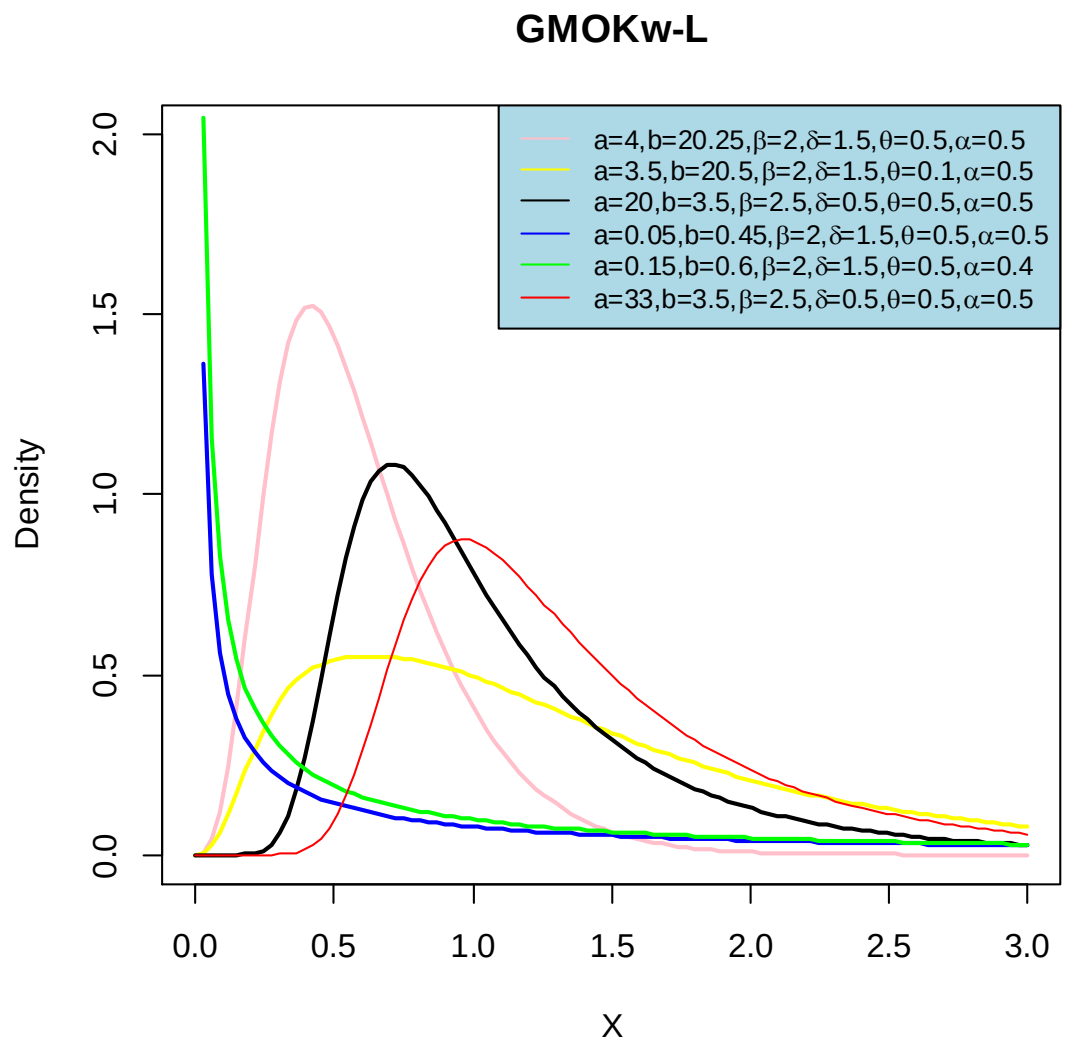

(c)

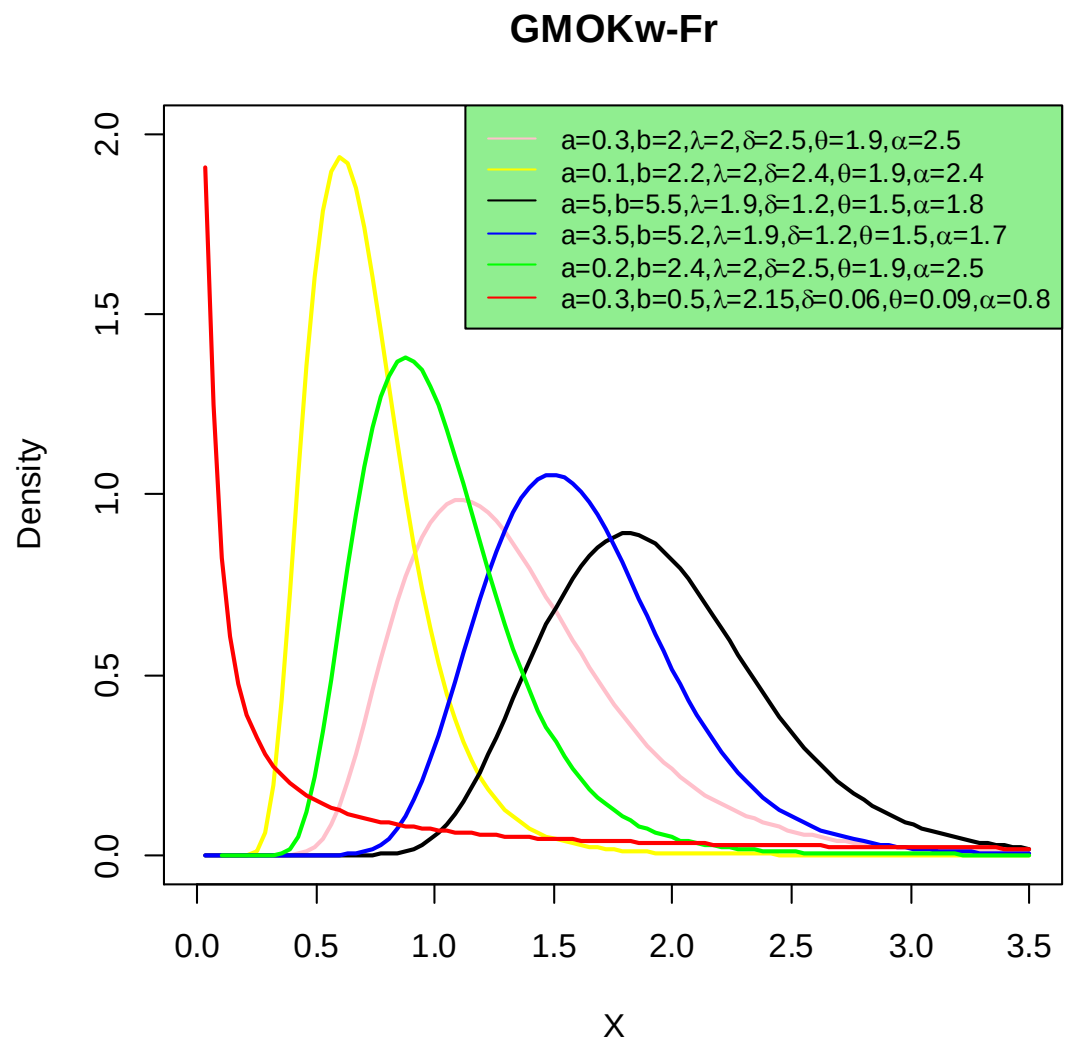

(d)

**Fig 1:** Density plots of (a) $GMOKw-E$, (b) $GMOKw-W$, (c) $GMOKw-L$ and (d) $GMOKw-Fr$ distributions.

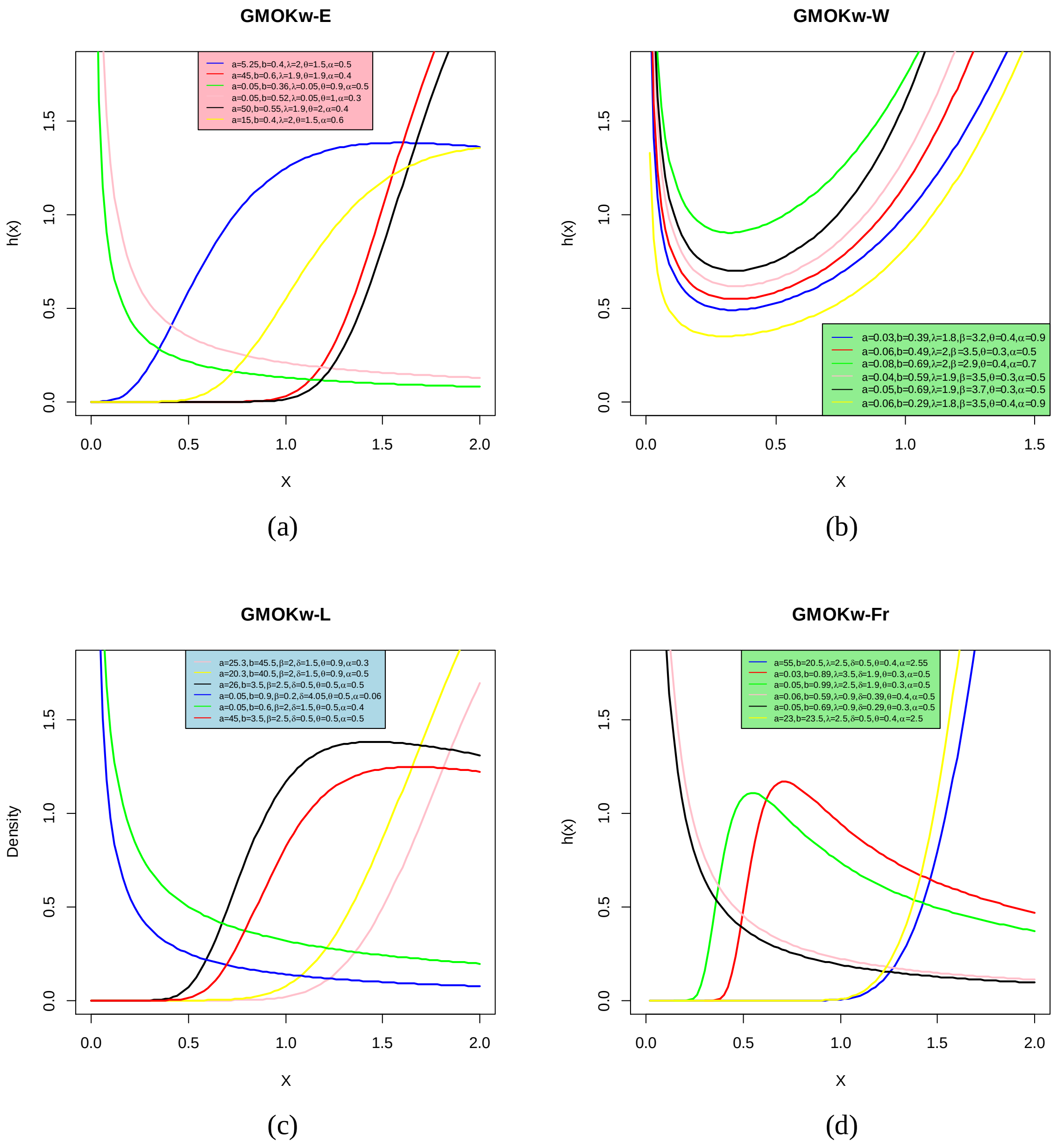

**Fig 2:** Hazard rate function plots of (a) $GMOKw-E$, (b) $GMOKw-W$, (c) $GMOKw-L$ and (d) $GMOKw-Fr$ distributions.

From the plots in figure 1 and 2 it can be seen that the family is very flexible and can offer many different types of shapes. It offers IFR, DFR even bath tub shaped hazard rate.

## 4. Some special $GMOKw-G$ distribution

Some special cases of the $GMOKw-G$ family of distributions are presented in this section.

### 4.1 The $GMOKw-$exponential ($GMOKw-E$) distribution

Let the base line distribution be exponential with parameter $\lambda>0$, $g(t:\lambda)=\lambda e^{-\lambda t}$, $t>0$ and $G(t:\lambda)=1-e^{-\lambda t}$, $t>0$ then for the $GMOKw-E$ model we get

pdf: $$f^{GMOKwE}(t)=\frac{\theta\alpha^{\theta}ab\lambda e^{-\lambda t}[1-e^{-\lambda t}]^{a-1}[1-[1-e^{-\lambda t}]^{a}]^{b\theta-1}}{[1-\bar{\alpha}[1-[1-e^{-\lambda t}]^{a}]^{b}]^{\theta+1}}$$

cdf: $$F^{GMOKwE}(t)=1-\left[\frac{\alpha[1-\{1-e^{-\lambda t}\}^{a}]^{b}}{1-\bar{\alpha}[1-\{1-e^{-\lambda t}\}^{a}]^{b}}\right]^{\theta}$$

sf: $$\bar{F}^{GMOKwE}(t)=\left[\frac{\alpha[1-\{1-e^{-\lambda t}\}^{a}]^{b}}{1-\bar{\alpha}[1-\{1-e^{-\lambda t}\}^{a}]^{b}}\right]^{\theta}$$

hrf : $$h^{GMOKwE}(t)=\frac{\theta ab\lambda e^{-\lambda t}[1-e^{-\lambda t}]^{a-1}[1-\{1-e^{-\lambda t}\}^{a}]^{-1}}{1-\bar{\alpha}[1-\{1-e^{-\lambda t}\}^{a}]^{b}}$$

rhrf: $$r^{GMOKwE}(t)=\frac{\theta\alpha^{\theta}ab\lambda e^{-\lambda t}[1-e^{-\lambda t}]^{a-1}[1-\{1-e^{-\lambda t}\}^{a}]^{b\theta-1}}{[1-\bar{\alpha}[1-\{1-e^{-\lambda t}\}^{a}]^{b}]^{\theta+1}-\alpha^{\theta}[1-\{1-e^{-\lambda t}\}^{a}]^{b\theta}[1-\bar{\alpha}[1-\{1-e^{-\lambda t}\}^{a}]^{b}]}$$

chrf: $$H^{GMOKwE}(t)=-\theta\log\left[\frac{\alpha[1-\{1-e^{-\lambda t}\}^{a}]^{b}}{1-\bar{\alpha}[1-\{1-e^{-\lambda t}\}^{a}]^{b}}\right]$$

### 4.2 The $GMOKw-$Lomax ($GMOKw-L$) distribution

Considering the Lomax distribution (Ghitany *et al.* 2007) with respective pdf and cdf $g(t:\beta,\delta)=(\beta/\delta)[1+(t/\delta)]^{-(\beta+1)}$, $t>0$, and $G(t:\beta,\delta)=1-[1+(t/\delta)]^{-\beta}$, $\beta>0$, $\delta>0$, the $GMOKw-L$ distribution is obtained with

pdf: $$f^{GMOKwL}(t)=\frac{\theta\alpha^{\theta}ab(\beta/\delta)[1+(t/\delta)]^{-(\beta+1)}[1-[1+(t/\delta)]^{-\beta}]^{a-1}[1-[1-[1+(t/\delta)]^{-\beta}]^{a}]^{b\theta-1}}{[1-\bar{\alpha}[1-[1-[1+(t/\delta)]^{-\beta}]^{a}]^{b}]^{\theta+1}}$$

cdf: $$F^{GMOKwL}(t)=1-\left[\frac{\alpha[1-\{1-[1+(t/\delta)]^{-\beta}\}^{a}]^{b}}{1-\bar{\alpha}[1-\{1-[1+(t/\delta)]^{-\beta}\}^{a}]^{b}}\right]^{\theta}$$

sf: $$\bar{F}^{GMOKwL}(t)=\left[\frac{\alpha[1-\{1-[1+(t/\delta)]^{-\beta}\}^{a}]^{b}}{1-\bar{\alpha}[1-\{1-[1+(t/\delta)]^{-\beta}\}^{a}]^{b}}\right]^{\theta}$$

hrf: $$h^{GMOKwL}(t)=\frac{\theta ab(\beta/\delta)[1+(t/\delta)]^{-(\beta+1)}[1-[1+(t/\delta)]^{-\beta}]^{a-1}[1-\{1-[1+(t/\delta)]^{-\beta}\}^{a}]^{-1}}{1-\bar{\alpha}[1-\{1-[1+(t/\delta)]^{-\beta}\}^{a}]^{b}}$$

rhrf : $r^{GMOKwL}(t)=$

$$\frac{\theta\alpha^{\theta}ab(\beta/\delta)[1+(t/\delta)]^{-(\beta+1)}\{1-[1+(t/\delta)]^{-\beta}\}^{a-1}}{[1-\bar{\alpha}[1-\{1-[1+(t/\delta)]^{-\beta}\}^{a}]^{b}]^{\theta+1}-\alpha^{\theta}[[1-\{1-[1+(t/\delta)]^{-\beta}\}^{a}]^{b\theta}}$$

$$\frac{[1-\{1-[1+(t/\delta)]^{-\beta}\}^{a}]^{b\theta-1}}{[1-\bar{\alpha}[1-\{1-[1+(t/\delta)]^{-\beta}\}^{a}]^{b}]]}$$

chrf: $$H^{GMOKwL}(t) = -\theta \log\left[\frac{\alpha[1-\{1-[1+(t/\delta)]^{-\beta}\}^{a}]^{b}}{1-\bar{\alpha}[1-\{1-[1+(t/\delta)]^{-\beta}\}^{a}]^{b}}\right]$$

**4.3 The $GMOKw-$Weibull ($GMOKw-W$) distribution**

Considering the Weibull distribution (Ghitany *et al.* 2005, Zhang and Xie 2007) with parameters $\lambda > 0$ and $\beta > 0$ having pdf and cdf $g(t) = \lambda\beta\, t^{\beta-1} e^{-\lambda t^{\beta}}$ and $G(t) = 1-e^{-\lambda t^{\beta}}$ respectively we get the $GMOKw-W$ distribution with

pdf: $$f^{GMOKwW}(t) = \frac{\theta\alpha^{\theta} ab\lambda\beta t^{\beta-1} e^{-\lambda t^{\beta}}[1-e^{-\lambda t^{\beta}}]^{a-1}[1-\{1-e^{-\lambda t^{\beta}}\}^{a}]^{b\theta-1}}{[1-\bar{\alpha}[1-\{1-e^{-\lambda t^{\beta}}\}^{a}]^{b}]^{\theta+1}}$$

cdf: $$F^{GMOKwW}(t) = 1-\left[\frac{\alpha[1-\{1-e^{-\lambda t^{\beta}}\}^{a}]^{b}}{1-\bar{\alpha}[1-\{1-e^{-\lambda t^{\beta}}\}^{a}]^{b}}\right]^{\theta}$$

sf: $$\bar{F}^{GMOKwW}(t) = \left[\frac{\alpha[1-\{1-e^{-\lambda t^{\beta}}\}^{a}]^{b}}{1-\bar{\alpha}[1-\{1-e^{-\lambda t^{\beta}}\}^{a}]^{b}}\right]^{\theta}$$

hrf: $$h^{GMOKwW}(t) = \frac{\theta ab\lambda\beta\, t^{\beta-1} e^{-\lambda t^{\beta}}[1-e^{-\lambda t^{\beta}}]^{a-1}[1-\{1-e^{-\lambda t^{\beta}}\}^{a}]^{-1}}{1-\bar{\alpha}[1-\{1-e^{-\lambda t^{\beta}}\}^{a}]^{b}}$$

rhrf : $r^{GMOKwW}(t) =$

$$\frac{\theta\alpha^{\theta} ab\lambda\beta\, t^{\beta-1} e^{-\lambda t^{\beta}}[1-e^{-\lambda t^{\beta}}]^{a-1}[1-\{1-e^{-\lambda t^{\beta}}\}^{a}]^{b\theta-1}}{[1-\bar{\alpha}[1-\{1-e^{-\lambda t^{\beta}}\}^{a}]^{b}]^{\theta+1} - \alpha^{\theta}[1-\{1-e^{-\lambda t^{\beta}}\}^{a}]^{b\theta}[1-\bar{\alpha}[1-\{1-e^{-\lambda t^{\beta}}\}^{a}]^{b}]}$$

chrf: $$H^{GMOKwW}(t) = -\theta \log\left[\frac{\alpha[1-\{1-e^{-\lambda t^{\beta}}\}^{a}]^{b}}{1-\bar{\alpha}[1-\{1-e^{-\lambda t^{\beta}}\}^{a}]^{b}}\right]$$

**4.4 The $GMOKw-$Frechet ($GMOKw-Fr$) distribution**

Suppose the base line distribution is the Frechet distribution (Krishna *et al.*, 2013) with pdf and cdf given by $g(t) = \lambda\delta^{\lambda} t^{-(\lambda+1)} e^{-(\delta/t)^{\lambda}}$ and $G(t) = e^{-(\delta/t)^{\lambda}}, t > 0$ respectively, and then the corresponding pdf of $GMOKw-Fr$ distribution is obtained with

pdf: $$f^{GMOKwFr}(t) = \frac{\theta\alpha^{\theta} ab\lambda\delta^{\lambda} t^{-(\lambda+1)} e^{-(\delta/t)^{\lambda}}[e^{-(\delta/t)^{\lambda}}]^{a-1}[1-\{e^{-(\delta/t)^{\lambda}}\}^{a}]^{b\theta-1}}{[1-\bar{\alpha}[1-\{e^{-(\delta/t)^{\lambda}}\}^{a}]^{b}]^{\theta+1}}$$

cdf: $$F^{GMOKwFr}(t) = 1-\left[\frac{\alpha[1-\{e^{-(\delta/t)^{\lambda}}\}^{a}]^{b}}{1-\bar{\alpha}[1-\{e^{-(\delta/t)^{\lambda}}\}^{a}]^{b}}\right]^{\theta}$$

sf: $\bar{F}^{GMOKwFr}(t)=\left[\dfrac{\alpha[1-\{e^{-(\delta/t)^{\lambda}}\}^{a}]^{b}}{1-\bar{\alpha}[1-\{e^{-(\delta/t)^{\lambda}}\}^{a}]^{b}}\right]^{\theta}$

hrf: $h^{GMOKwFr}(t)=\dfrac{\theta ab\lambda\delta^{\lambda}t^{-(\lambda+1)}e^{-(\delta/t)^{\lambda}}[e^{-(\delta/t)^{\lambda}}]^{a-1}[1-\{e^{-(\delta/t)^{\lambda}}\}^{a}]^{-1}}{1-\bar{\alpha}[1-\{e^{-(\delta/t)^{\lambda}}\}^{a}]^{b}}$

rhrf: $r^{GMOKwFr}(t)=:\ \dfrac{\theta\alpha^{\theta}ab\lambda\delta^{\lambda}t^{-(\lambda+1)}e^{-(\delta/t)^{\lambda}}[e^{-(\delta/t)^{\lambda}}]^{a-1}[1-\{e^{-(\delta/t)^{\lambda}}\}^{a}]^{b\theta-1}}{[1-\bar{\alpha}[1-\{e^{-(\delta/t)^{\lambda}}\}^{a}]^{b}]^{\theta+1}-\alpha^{\theta}[1-\{e^{-(\delta/t)^{\lambda}}\}^{a}]^{b\theta}[1-\bar{\alpha}[1-\{e^{-(\delta/t)^{\lambda}}\}^{a}]^{b}]}$

chrf: $H^{GMOKwFr}(t)=-\theta\log\left[\dfrac{\alpha[1-\{e^{-(\delta/t)^{\lambda}}\}^{a}]^{b}}{1-\bar{\alpha}[1-\{e^{-(\delta/t)^{\lambda}}\}^{a}]^{b}}\right]$

### 4.5 The $GMOKw-$Gompertz ($GMOKw-Go$) distribution

Considering the Gompertz distribution (Gieser et al. 1998) with pdf and cdf $g(t)=\beta e^{\lambda t}e^{-\frac{\beta}{\lambda}(e^{\lambda t}-1)}$ and $G(t)=1-e^{-\frac{\beta}{\lambda}(e^{\lambda t}-1)}$, $\beta>0,\ \lambda>0,\ t>0$ respectively, we get the $GMOKw-Go$ distribution with pdf: $f^{GMOKwGo}(t)=\dfrac{\theta\alpha^{\theta}ab\beta e^{\lambda t}e^{-\frac{\beta}{\lambda}(e^{\lambda t}-1)}[1-e^{-\frac{\beta}{\lambda}(e^{\lambda t}-1)}]^{a-1}[1-\{1-e^{-\frac{\beta}{\lambda}(e^{\lambda t}-1)}\}^{a}]^{b\theta-1}}{[1-\bar{\alpha}[1-\{1-e^{-\frac{\beta}{\lambda}(e^{\lambda t}-1)}\}^{a}]^{b}]^{\theta+1}}$ cdf: $F^{GMOKwGo}(t)=1-\left[\dfrac{\alpha[1-\{1-e^{-\frac{\beta}{\lambda}(e^{\lambda t}-1)}\}^{a}]^{b}}{1-\bar{\alpha}[1-\{1-e^{-\frac{\beta}{\lambda}(e^{\lambda t}-1)}\}^{a}]^{b}}\right]^{\theta}$ sf: $\bar{F}^{GMOKwGo}(t)=\left[\dfrac{\alpha[1-\{1-e^{-\frac{\beta}{\lambda}(e^{\lambda t}-1)}\}^{a}]^{b}}{1-\bar{\alpha}[1-\{1-e^{-\frac{\beta}{\lambda}(e^{\lambda t}-1)}\}^{a}]^{b}}\right]^{\theta}$ hrf: $h^{GMOKwGo}(t)=\dfrac{\theta ab\beta e^{\lambda t}e^{-\frac{\beta}{\lambda}(e^{\lambda t}-1)}[1-e^{-\frac{\beta}{\lambda}(e^{\lambda t}-1)}]^{a-1}[1-\{1-e^{-\frac{\beta}{\lambda}(e^{\lambda t}-1)}\}^{a}]^{-1}}{1-\bar{\alpha}[1-\{1-e^{-\frac{\beta}{\lambda}(e^{\lambda t}-1)}\}^{a}]^{b}}$ rhrf: $r^{GMOKwGo}(t)=\dfrac{\theta\alpha^{\theta}ab\beta e^{\lambda t}e^{-\frac{\beta}{\lambda}(e^{\lambda t}-1)}[1-e^{-\frac{\beta}{\lambda}(e^{\lambda t}-1)}]^{a-1}}{[1-\bar{\alpha}[1-\{1-e^{-\frac{\beta}{\lambda}(e^{\lambda t}-1)}\}^{a}]^{b}]^{\theta+1}-\alpha^{\theta}[[1-\{1-e^{-\frac{\beta}{\lambda}(e^{\lambda t}-1)}\}^{a}]^{b\theta}}$

$$\times\frac{[1-\{1-e^{-\frac{\beta}{\lambda}(e^{\lambda t}-1)}\}^{a}]^{b\theta-1}}{[1-\bar{\alpha}[1-\{1-e^{-\frac{\beta}{\lambda}(e^{\lambda t}-1)}\}^{a}]^{b}]]}$$

chrf: $H^{GMOKwGo}(t)=-\theta\log\left[\dfrac{\alpha[1-\{1-e^{-\frac{\beta}{\lambda}(e^{\lambda t}-1)}\}^{a}]^{b}}{1-\bar{\alpha}[1-\{1-e^{-\frac{\beta}{\lambda}(e^{\lambda t}-1)}\}^{a}]^{b}}\right]$

**4.6 The $GMOKw-$Extended Weibull ($GMOKw-EW$) distribution**

The extended Weibull (*EW*) distributions of Gurvich *et al.* (1997) is defined by cdf $G(t:\delta,\xi)=1-\exp[-\delta\, E(t:\xi)]$, $\quad t\in D\subseteq R_{+},\delta>0$ where $E(t:\xi)$ is a non-negative monotonically increasing function which depends on the parameter vector $\xi$**.** The corresponding pdf is $g(t:\delta,\xi)=\delta\exp[-\delta\, E(t:\xi)]\, e(t:\xi)$ where $e(t:\xi)$ is of course the derivative of $E(t:\xi)$. Many important distributions can be obtained by choosing different expressions for $E(t:\xi)$.

The $GMOKw-EW$ is derived by considering *EW* as the base line distribution with

pdf: $$f^{GMOKwEW}(t)=\frac{\theta\alpha^{\theta}ab\delta\exp[-\delta\, E(t:\xi)]e(t:\xi)[1-\exp\{-\delta\, E(t:\xi)\}]^{a-1}}{[1-\overline{\alpha}[1-\{1-\exp[-\delta\, E(t:\xi)]\}^{a}]^{b}]^{\theta+1}}$$
$$\times[1-\{1-\exp[-\delta\, E(t:\xi)]\}^{a}]^{b\theta-1}$$

cdf: $$F^{GMOKwEW}(t)=1-\left[\frac{\alpha[1-\{1-\exp[-\delta\, E(t:\xi)]\}^{a}]^{b}}{1-\overline{\alpha}[1-\{1-\exp[-\delta\, E(t:\xi)]\}^{a}]^{b}}\right]^{\theta}$$

sf: $$\overline{F}^{GMOKwEW}(t)=\left[\frac{\alpha[1-\{1-\exp[-\delta\, E(t:\xi)]\}^{a}]^{b}}{1-\overline{\alpha}[1-\{1-\exp[-\delta\, E(t:\xi)]\}^{a}]^{b}}\right]^{\theta}$$

hrf: $$h^{GMOKwEW}(t)=\frac{\theta ab\delta\exp[-\delta\, E(t:\xi)]e(t:\xi)\{1-\exp[-\delta\, E(t:\xi)]\}^{a-1}}{1-\overline{\alpha}[1-\{1-\exp[-\delta\, E(t:\xi)]\}^{a}]^{b}}$$
$$[1-\{1-\exp[-\delta\, E(t:\xi)]\}^{a}]^{-1}$$

rhrf: $r^{GMOKwEW}(t)=$

$$\frac{\theta\alpha^{\theta}ab\delta\exp[-\delta\, E(t:\xi)]e(t:\xi)\{1-\exp[-\delta\, E(t:\xi)]\}^{a-1}}{[1-\overline{\alpha}[1-\{1-\exp[-\delta\, E(t:\xi)]\}^{a}]^{b}]^{\theta+1}-\alpha^{\theta}[[1-\{1-\exp[-\delta\, E(t:\xi)]\}^{a}]^{b\theta}}$$
$$\frac{[1-\{1-\exp[-\delta\, E(t:\xi)]\}^{a}]^{b\theta-1}}{[1-\overline{\alpha}[1-\{1-\exp[-\delta\, E(t:\xi)]\}^{a}]^{b}]]}$$

chrf: $$H^{GMOKwEW}(t)=-\theta\log\left[\frac{\alpha[1-\{1-\exp[-\delta\, E(t:\xi)]\}^{a}]^{b}}{1-\overline{\alpha}[1-\{1-\exp[-\delta\, E(t:\xi)]\}^{a}]^{b}}\right].$$

**4.7 The $GMOKw-$Extended Modified Weibull ($GMOKw-EMW$) distribution**

The modified Weibull (*MW*) distribution of Sarhan and Zaindin (2013) has cdf and pdf is given by

$G(t;\sigma,\beta,\gamma)=1-\exp[-\sigma t-\beta t^{\gamma}]$ **,** $t>0,\gamma>0,\sigma,\beta\geq 0,\sigma+\beta>0$ and

$g(t;\sigma,\beta,\gamma)=(\sigma+\beta\gamma t^{\gamma-1})\exp[-\sigma t-\beta t^{\gamma}]$ respectively.

The corresponding $GMOKw-EMW$ is obtained with pdf: $f^{GMOKwEMW}(t) = \theta\alpha^{\theta} ab(\sigma+\beta\gamma t^{\gamma-1})\exp[-\sigma t-\beta t^{\gamma}]\{1-\exp[-\sigma t-\beta t^{\gamma}]\}^{a-1}$

$$\frac{[1-\{1-\exp[-\sigma t-\beta t^{\gamma}]\}^{a}]^{b\theta-1}}{[1-\bar{\alpha}[1-\{1-\exp[-\sigma t-\beta t^{\gamma}]\}^{a}]^{b}]^{\theta+1}}$$

cdf: $F^{GMOKwEMW}(t) = 1-\left[\frac{\alpha[1-\{1-\exp[-\sigma t-\beta t^{\gamma}]\}^{a}]^{b}}{1-\bar{\alpha}[1-\{1-\exp[-\sigma t-\beta t^{\gamma}]\}^{a}]^{b}}\right]^{\theta}$

sf: $\bar{F}^{GMOKwEMW}(t) = \left[\frac{\alpha[1-\{1-\exp[-\sigma t-\beta t^{\gamma}]\}^{a}]^{b}}{1-\bar{\alpha}[1-\{1-\exp[-\sigma t-\beta t^{\gamma}]\}^{a}]^{b}}\right]^{\theta}$

hrf: $h^{GMOKwEMW}(t) = \frac{\theta ab(\sigma+\beta\gamma t^{\gamma-1})\exp[-\sigma t-\beta t^{\gamma}]\{1-\exp[-\sigma t-\beta t^{\gamma}]\}^{a-1}}{1-\bar{\alpha}[1-\{1-\exp[-\sigma t-\beta t^{\gamma}]\}^{a}]^{b}}$

$$\times[1-\{1-\exp[-\sigma t-\beta t^{\gamma}]\}^{a}]^{-1}$$

rhrf : $r^{GMOKwEMW}(t) =$

$$\frac{\theta\alpha^{\theta} ab(\sigma+\beta\gamma t^{\gamma-1})\exp[-\sigma t-\beta t^{\gamma}][1-\exp[-\sigma t-\beta t^{\gamma}]]^{a-1}}{[1-\bar{\alpha}[1-\{1-\exp[-\sigma t-\beta t^{\gamma}]\}^{a}]^{b}]^{\theta+1}-\alpha^{\theta}[[1-\{1-\exp[-\sigma t-\beta t^{\gamma}]\}^{a}]^{b\theta}}$$

$$\frac{[1-\{1-\exp[-\sigma t-\beta t^{\gamma}]]^{a}]^{b\theta-1}}{[1-\bar{\alpha}[1-\{1-\exp[-\sigma t-\beta t^{\gamma}]\}^{a}]^{b}]]}$$

chrf: $H^{GMOKwEMW}(t) = -\theta\log\left[\frac{\alpha[1-\{1-\exp[-\sigma t-\beta t^{\gamma}]\}^{a}]^{b}}{1-\bar{\alpha}[1-\{1-\exp[-\sigma t-\beta t^{\gamma}]\}^{a}]^{b}}\right]$

**4.8 The $GMOKw-$Extended Exponentiated Pareto ($GMOKw-EEP$) distribution**

The pdf and cdf of the exponentiated Pareto distribution of Nadarajah (2005), are given respectively by $g(t)=\gamma k\theta^{k}t^{-(k+1)}[1-(\theta/t)^{k}]^{\gamma-1}$ and $G(t)=[1-(\theta/t)^{k}]^{\gamma}$, $x>\theta$ and $\theta, k, \gamma>0$. Thus the $GMOKw-EEP$ distribution is defined with pdf: $f^{GMOKwEEP}(t) = \frac{\theta\alpha^{\theta} ab\gamma k\theta^{k}t^{-(k+1)}[1-(\theta/t)^{k}]^{\gamma-1}\{[1-(\theta/t)^{k}]^{\gamma}\}^{a-1}[1-\{[1-(\theta/t)^{k}]^{\gamma}\}^{a}]^{b\theta-1}}{[1-\bar{\alpha}[1-\{[1-(\theta/t)^{k}]^{\gamma}\}^{a}]^{b}]^{\theta+1}}$

cdf: $F^{GMOKwEEP}(t) = 1-\left[\frac{\alpha[1-\{[1-(\theta/t)^{k}]^{\gamma}\}^{a}]^{b}}{1-\bar{\alpha}[1-\{[1-(\theta/t)^{k}]^{\gamma}\}^{a}]^{b}}\right]^{\theta}$

sf: $\bar{F}^{GMOKwEEP}(t) = \left[\frac{\alpha[1-\{[1-(\theta/t)^{k}]^{\gamma}\}^{a}]^{b}}{1-\bar{\alpha}[1-\{[1-(\theta/t)^{k}]^{\gamma}\}^{a}]^{b}}\right]^{\theta}$

hrf: $h^{GMOKwEEP}(t) = \frac{\theta ab\gamma k\theta^{k}t^{-(k+1)}[1-(\theta/t)^{k}]^{\gamma-1}\{[1-(\theta/t)^{k}]^{\gamma}\}^{a-1}[1-\{[1-(\theta/t)^{k}]^{\gamma}\}^{a}]^{-1}}{1-\bar{\alpha}[1-\{[1-(\theta/t)^{k}]^{\gamma}\}^{a}]^{b}}$

rhrf: $$r^{GMOKwEEP}(t)=\frac{\theta\alpha^{\theta}ab\gamma k\theta^{k}t^{-(k+1)}[1-(\theta/t)^{k}]^{\gamma-1}\{[1-(\theta/t)^{k}]^{\gamma}\}^{a-1}}{[1-\overline{\alpha}[1-\{[1-(\theta/t)^{k}]^{\gamma}\}^{a}]^{b}]^{\theta+1}-\alpha^{\theta}[[1-\{[1-(\theta/t)^{k}]^{\gamma}\}^{a}]^{b\theta}}$$

$$\frac{[1-\{[1-(\theta/t)^{k}]^{\gamma}\}^{a}]^{b\theta-1}}{[1-\overline{\alpha}[1-\{[1-(\theta/t)^{k}]^{\gamma}\}^{a}]^{b}]]}$$

chrf: $$H^{GMOKwEEP}(t)=-\theta\log\left[\frac{\alpha[1-\{[1-(\theta/t)^{k}]^{\gamma}\}^{a}]^{b}}{1-\overline{\alpha}[1-\{[1-(\theta/t)^{k}]^{\gamma}\}^{a}]^{b}}\right]$$

### 4.9 The $GMOKw-$Extended Power ($GMOKw-EP$) distribution

The cdf and pdf of the extended power distribution are respectively given by $G(t)=(\theta t)^{k}$ and $g(t)=k\theta^{k}t^{k-1}$, $t\in(0,1/\theta)$ and $\theta>0$.

The corresponding $GMOKw-EP$ distribution is then given by

pdf: $$f^{GMOKwEP}(t)=\frac{\theta\alpha^{\theta}abk\theta^{k}t^{k-1}\{(\theta t)^{k}\}^{a-1}[1-\{(\theta t)^{k}\}^{a}]^{b\theta-1}}{[1-\overline{\alpha}[1-\{(\theta t)^{k}\}^{a}]^{b}]^{\theta+1}}$$

cdf: $$F^{GMOKwEP}(t)=1-\left[\frac{\alpha[1-\{(\theta t)^{k}\}^{a}]^{b}}{1-\overline{\alpha}[1-\{(\theta t)^{k}\}^{a}]^{b}}\right]^{\theta}$$

sf: $$\overline{F}^{GMOKwEP}(t)=\left[\frac{\alpha[1-\{(\theta t)^{k}\}^{a}]^{b}}{1-\overline{\alpha}[1-\{(\theta t)^{k}\}^{a}]^{b}}\right]^{\theta}$$

hrf: $$h^{GMOKwEP}(t)=\frac{\theta abk\theta^{k}t^{k-1}\{(\theta t)^{k}\}^{a-1}[1-\{(\theta t)^{k}\}^{a}]^{-1}}{1-\overline{\alpha}[1-\{(\theta t)^{k}\}^{a}]^{b}}$$

rhrf: $$r^{GMOKwEP}(t)=\frac{\theta\alpha^{\theta}abk\theta^{k}t^{k-1}\{(\theta t)^{k}\}^{a-1}[1-\{(\theta t)^{k}\}^{a}]^{b\theta-1}}{[1-\overline{\alpha}[1-\{(\theta t)^{k}\}^{a}]^{b}]^{\theta+1}-\alpha^{\theta}[1-\{(\theta t)^{k}\}^{a}]^{b\theta}[1-\overline{\alpha}[1-\{(\theta t)^{k}\}^{a}]^{b}]}$$

chrf: $$H^{GMOKwEP}(t)=-\theta\log\left[\frac{\alpha[1-\{(\theta t)^{k}\}^{a}]^{b}}{1-\overline{\alpha}[1-\{(\theta t)^{k}\}^{a}]^{b}}\right]$$

## 5. General results for the Generalized Marshall-Olkin Kumaraswamy-*G* ($GMOKw-G$) family of distributions

In this section we derive some general results for the proposed $GMOKw-G$ family in a similar manner as done in Barreto-Souza *et al.* (2013), Cordeiro *et al.* (2014), and Alizadeh *et al.* (2015) among others.

### 5.1 Expansions

We know that $(1-z)^{-k}=\sum_{j=0}^{\infty}\frac{\Gamma(k+j)}{\Gamma(k)\,j!}z^{j}$ , $|z|<1$ and $k>0$ (9)

where $\Gamma(.)$ is the gamma function.

If $\alpha\in(0,1)$ using (9) in (6), we obtain

$$f^{GMOKwG}(t)=\theta\alpha^{\theta}\,a\,b\,g(t)G(t)^{a-1}[1-G(t)^{a}]^{b\theta-1}[1-\bar{\alpha}[1-G(t)^{a}]^{b}]^{-(\theta+1)}$$

$$=\theta\alpha^{\theta}\,a\,b\,g(t)G(t)^{a-1}[1-G(t)^{a}]^{b\theta-1}\sum_{j=0}^{\infty}\frac{\Gamma(j+\theta+1)}{\Gamma(\theta+1)\,j!}\{\bar{\alpha}[1-G(t)^{a}]^{b}\}^{j}$$

$$=\theta\alpha^{\theta}\,a\,b\,g(t)G(t)^{a-1}[1-G(t)^{a}]^{b\theta-1}\sum_{j=0}^{\infty}\frac{(j+\theta)!}{\theta!\,j!}(1-\alpha)^{j}\{[1-G(t)^{a}]^{b}\}^{j}$$

$$=a\,b\,g(t)G(t)^{a-1}[1-G(t)^{a}]^{b\theta-1}\sum_{j=0}^{\infty}A_{j}\{[1-G(t)^{a}]^{b}\}^{j}$$

$$=[1-G(t)^{a}]^{b(\theta-1)}f^{KwG}(t;a,b)\sum_{j=0}^{\infty}A_{j}[\bar{F}^{KwG}(t;a,b)]^{j}$$

$$=f^{KwG}(t;a,b)\sum_{j=0}^{\infty}A_{j}[\bar{F}^{KwG}(t;a,b)]^{j+\theta-1} \qquad (10)$$

$$=\sum_{j=0}^{\infty}\frac{A_{j}}{j+\theta}\frac{d}{dt}[\bar{F}^{KwG}(t;a,b)]^{j+\theta}$$

$$=\sum_{j=0}^{\infty}A'_{j}\frac{d}{dt}[\bar{F}^{KwG}(t;a,b)]^{j+\theta} \qquad (11)$$

$$=\sum_{j=0}^{\infty}A'_{j}\frac{d}{dt}[\bar{F}^{KwG}(t;a,b(j+\theta))]$$

$$=\sum_{j=0}^{\infty}A'_{j}\,f^{KwG}(t;a,b(j+\theta)) \qquad (12)$$

Where $A'_{j}=A'_{j}(\alpha)=-\binom{j+\theta-1}{j}(1-\alpha)^{j}\alpha^{\theta}$, $\quad A_{j}=A_{j}(\alpha)=(j+\theta)A'_{j}$

From this formulation the $\text{GMOKw-G}(\alpha,\theta,a,b)$ can be called a non central $Kw-G$ (Cordeiro and de Castro, 2011) distribution. Similarly an expansion for the survival function of $GMOKw-G$ [for $\alpha\in(0,1)$ ] can be derives as

$$\bar{F}^{GMOKwG}(t)=\left[\frac{\alpha[1-G(t)^{a}]^{b}}{1-\bar{\alpha}[1-G(t)^{a}]^{b}}\right]^{\theta}$$

$$=\{\alpha[1-G(t)^{a}]^{b}\}^{\theta}\{1-\bar{\alpha}[1-G(t)^{a}]^{b}\}^{-\theta}$$

$$=\{\alpha[1-G(t)^{a}]^{b}\}^{\theta}\sum_{j=0}^{\infty}\frac{\Gamma(j+\theta)}{\Gamma(\theta)\,j!}\{\bar{\alpha}[1-G(t)^{a}]^{b}\}^{j}$$

$$=\alpha^{\theta}\{\bar{F}^{KwG}(t;a,b)\}^{\theta}\sum_{j=0}^{\infty}\frac{(j+\theta-1)!}{(\theta-1)!\,j!}(1-\alpha)^{j}\{\bar{F}^{KwG}(t;a,b)\}^{j}$$

$$= -\sum_{j=0}^{\infty} A'_j [\bar{F}^{KwG}(t;a,b)]^{j+\theta} = -\sum_{j=0}^{\infty} A'_j \bar{F}^{KwG}(t;a,b(j+\theta)) \quad (13)$$

Alternatively, we can expand the pdf as

$$f^{GMOKwG}(t) = \theta\alpha^{\theta} a b g(t) G(t)^{a-1}[1-G(t)^a]^{b\theta-1} \sum_{j=0}^{\infty} \frac{(j+\theta)!}{\theta!\, j!} (1-\alpha)^j [1-[1-\{[1-G(t)^a]^b\}^j]]$$

$$= \theta\alpha^{\theta} a b g(t) G(t)^{a-1}[1-G(t)^a]^{b\theta-1} \sum_{j=0}^{\infty} \frac{(j+\theta)!}{\theta!\, j!} (1-\alpha)^j \sum_{k=0}^{j} \binom{j}{k} [-\{1-[1-G(t)^a]^b\}]^{j-k}$$

$$= \theta\alpha^{\theta} a b g(t) G(t)^{a-1}[1-G(t)^a]^{b\theta-1} \sum_{j=0}^{\infty} \frac{(j+\theta)!}{\theta!\, j!} (1-\alpha)^j (-1)^{j-k} \sum_{k=0}^{j} \binom{j}{k} [1-[1-G(t)^a]^b]^{j-k}$$

$$= a b g(t) G(t)^{a-1}[1-G(t)^a]^{b\theta-1} \sum_{j=0}^{\infty} \sum_{k=0}^{j} \mathrm{B}_{j,k} [1-[1-G(t)^a]^b]^{j-k}$$

$$= f^{KwG}(t;a,b\theta) \sum_{j=0}^{\infty} \sum_{k=0}^{j} \mathrm{B}_{j,k} [F^{KwG}(t;a,b)]^{j-k} \quad (14)$$

Where $\mathrm{B}_{j,k} = \mathrm{B}_{j,k}(\alpha) = \frac{\theta\alpha^{\theta}(j+\theta)!}{\theta!\, j!}(1-\alpha)^j (-1)^{j-k} \binom{j}{k} = (-1)^{j-k-1} \binom{j}{k} \mathrm{A}_j$

Another expansion of the density function in (6) can be obtained by expressing the pdf as

$$f^{GMOKwG}(t) = \frac{\theta a b g(t) G(t)^{a-1}[1-G(t)^a]^{b\theta-1}}{\alpha\left[1-\frac{(\alpha-1)[1-[1-G(t)^a]^b]}{\alpha}\right]^{\theta+1}}$$

$$= \frac{\theta a b g(t) G(t)^{a-1}[1-G(t)^a]^{b\theta-1}}{\alpha}\left[1-\frac{(\alpha-1)[1-[1-G(t)^a]^b]}{\alpha}\right]^{-(\theta+1)}$$

Hence for $\alpha > 1$ using (9) we get

$$= \frac{\theta a b g(t) G(t)^{a-1}[1-G(t)^a]^{b\theta-1}}{\alpha^{j+1}} \sum_{j=0}^{\infty} \frac{\Gamma(j+\theta+1)}{\Gamma(\theta+1)\, j!} (\alpha-1)^j [1-[1-G(t)^a]^b]^j$$

$$= \frac{\theta a b g(t) G(t)^{a-1}[1-G(t)^a]^{b\theta-1}}{\alpha^{j+1}} \sum_{j=0}^{\infty} \frac{(j+\theta)!}{\theta!\, j!} (\alpha-1)^j [1-[1-G(t)^a]^b]^j$$

$$= a b g(t) G(t)^{a-1}[1-G(t)^a]^{b\theta-1} \sum_{j=0}^{\infty} C_j [1-[1-G(t)^a]^b]^j$$

$$= f^{KwG}(t;a,b\theta) \sum_{j=0}^{\infty} C_j \{F^{KwG}(t;a,b)\}^j \quad (15)$$

Where $C_j = C_j(\alpha) = \frac{\theta\,(j+\theta)!}{\theta!\, j!\, \alpha}\left(1-\frac{1}{\alpha}\right)^j = \frac{\theta}{\alpha}\binom{j+\theta-1}{j}\left(-\frac{\bar{\alpha}}{\alpha}\right)^j = \frac{(-1)^{j-1}}{\alpha^{\theta+j+1}} A_j$

For $\alpha > 1$, the survival function of $GMOKw-G$ can be expressed as

$$\overline{F}^{GMOKwG}(t)=\left[\frac{[1-G(t)^{a}]^{b}}{1-(\alpha-1)[1-[1-G(t)^{a}]^{b}]/\alpha}\right]^{\theta}$$

On using (9) we get

$$=\{[1-G(t)^{a}]^{b}\}^{\theta}\left[1-\frac{(\alpha-1)[1-[1-G(t)^{a}]^{b}]}{\alpha}\right]^{-\theta}$$

$$=\{[1-G(t)^{a}]^{b}\}^{\theta}\sum_{j=0}^{\infty}\frac{\Gamma(j+\theta)}{\Gamma(\theta)\ j!}\frac{(\alpha-1)^{j}}{\alpha^{j}}[1-[1-G(t)^{a}]^{b}]^{j}$$

$$=[\overline{F}^{KwG}(t;a,b)]^{\theta}\sum_{j=0}^{\infty}\binom{j+\theta-1}{j}(-1)^{j}\left(\frac{\overline{\alpha}}{\alpha}\right)^{j}[F^{KwG}(t;a,b)]^{j}$$

$$=[\overline{F}^{KwG}(t;a,b)]^{\theta}\sum_{j=0}^{\infty}C'_{j}[F^{KwG}(t;a,b)]^{j}$$

Where $C'_j=(\alpha/j+\theta)\ C_j$ .

**5.2 Order Statistics**

Suppose $T_1,T_2,...T_n$ is a random sample from any $GMOKw-G$ distribution. Let $T_{i:n}$ denote the *i*th order statistics. The pdf of $T_{i:n}$ can be expressed as

$$f_{i:n}(t)=\frac{n!}{(i-1)!\,(n-i)!}f^{GMOKwG}(t)[1-\overline{F}^{GMOKwG}(t)]^{i-1}\,\overline{F}^{GMOKwG}(t)^{n-i}$$

$$=\frac{n!}{(i-1)!\,(n-i)!}f^{GMOKwG}(t)\overline{F}^{GMOKwG}(t)^{n-i}\sum_{l=0}^{i-1}\binom{i-1}{l}[-\overline{F}^{GMOKwG}(t)]^{l}$$

$$=\frac{n!}{(i-1)!\,(n-i)!}f^{GMOKwG}(t)\sum_{l=0}^{i-1}(-1)^{l}\binom{i-1}{l}\overline{F}^{GMOKwG}(t)^{n+l-i}$$

Now using the general expansion of the $GMOKw-G$ distribution pdf and sf we get the pdf of the $i^{th}$ order statistics for of the $GMOKw-G$ for $\alpha\in(0,1)$ as

$$f_{i:n}(t)=\frac{n!}{(i-1)!\,(n-i)!}\left\{f^{KwG}(t;a,b)\sum_{j=0}^{\infty}A_j[\overline{F}^{KwG}(t;a,b)]^{j+\theta-1}\right\}$$

$$\sum_{l=0}^{i-1}(-1)^{l+1}\binom{i-1}{l}\left\{\sum_{p=0}^{\infty}A'_p[\overline{F}^{KwG}(t;a,b)]^{p+\theta(n+l-i)}\right\}$$

where $\mathrm{A}_j$ and $\mathrm{A}'_p$ are defined earlier

$$=\frac{n!}{(i-1)!\,(n-i)!}\sum_{l=o}^{i-1}\binom{i-1}{l}(-1)^{l+1}f^{KwG}(t;a,b)\sum_{j=0}^{\infty}\sum_{p=0}^{\infty}A_j\,A'_p[\overline{F}^{KwG}(t;a,b)]^{j+p+\theta(n+l-i+1)-1}$$

$$= f^{KwG}(t;a,b)\sum_{j,p=0}^{\infty} M_{j,p}[\bar{F}^{KwG}(t;a,b)]^{j+p+\theta(n+l-i+1)-1} \qquad (16)$$

$$= -\sum_{j,p=0}^{\infty}[M_{j,p}/(j+p+\theta(n+l-i+1))]\frac{d}{dt}[\bar{F}^{KwG}(t;a,b)]^{j+p+\theta(n+l-i+1)}$$

$$= \sum_{j,p=0}^{\infty} M'_{j,p}\frac{d}{dt}[\bar{F}^{KwG}(t;a,b)]^{j+p+\theta(n+l-i+1)}$$

$$= \sum_{j,p=0}^{\infty} M'_{j,p}\frac{d}{dt}[\bar{F}^{KwG}(t;a,b(j+p+\theta(n+l-i+1)))]$$

$$= \sum_{j,p=0}^{\infty} M'_{j,p}\, f^{KwG}(t;a,b(j+p+\theta(n+l-i+1)))] \qquad (17)$$

Where

$$M_{j,p} = nA_j A'_p \binom{n-1}{i-1}\sum_{l=0}^{i-1}\binom{i-1}{l}(-1)^{l+1} \text{ and } M'_{j,p} = -M_{j,p}/(j+p+\theta(n+l-i+1))$$

In particular for $\theta = 1$, (16) reduces to pdf of the $i^{th}$ order statistics for of the $\text{MOKw-G}(\alpha, a, b)$ given by

$$f_{i:n}(t) = f^{KwG}(t;a,b)\sum_{j,c=0}^{\infty} H_{j,c}[\bar{F}^{KwG}(t;a,b)]^{j+c+n+l-i}$$

Where $\quad H_{j,c} = n\kappa_j\,\eta_c\binom{n-1}{i-1}\sum_{l=0}^{i-1}\binom{i-1}{l}(-1)^{l+1}$

and $\quad \kappa_j = \kappa_j(\alpha) = (j+1)\alpha(1-\alpha)^j; \;\; \eta_c = \alpha(1-\alpha)^c$

Again using the general expansion of the $GMOKw-G$ distribution pdf and sf we get the pdf of the $i^{th}$ order statistics for of the $GMOKw-G$ for $\alpha > 1$ as

$$f_{i:n}(t) = \frac{n!}{(i-1)!\,(n-i)!}\left\{f^{KwG}(t;a,b\theta)\sum_{j=0}^{\infty} C_j\left\{F^{KwG}(t;a,b)\right\}^j\right\}\sum_{l=0}^{i-1}(-1)^l\binom{i-1}{l}$$
$$\left\{[\bar{F}^{KwG}(t;a,b)]^{\theta}\sum_{k=0}^{\infty} C'_k[F^{KwG}(t;a,b)]^k\right\}^{n+l-i}$$

$$= \frac{n!}{(i-1)!\,(n-i)!}\left\{f^{KwG}(t;a,b\theta)\sum_{j=0}^{\infty} C_j\left\{F^{KwG}(t;a,b)\right\}^j\right\}\sum_{l=0}^{i-1}(-1)^l\binom{i-1}{l}$$
$$[\bar{F}^{KwG}(t;a,b)]^{\theta(n+l-i)}\left[\sum_{k=0}^{\infty} C'_k[F^{KwG}(t;a,b)]^k\right]^{n+l-i}$$

Where $C_j$ and $C'_k$ defined in section 5.1

Now, $$\left[\sum_{k=0}^{\infty} C'_k [F^{KwG}(t;a,b)]^k\right]^{n+l-i} = \sum_{k=0}^{\infty} d_{n+l-i,k} [F^{KwG}(t;a,b)]^k \quad \text{(Nadarajah et. al 2015)}$$

Where $$d_{n+l-i,k} = \frac{1}{k C'_0} \sum_{h=1}^{k} [h(n+l-i-1)-k] C'_h d_{n+l-i,k-h}$$

Therefore the density function of the $i$th order statistics of $GMOKw-G$ distribution can be expressed as

$$f_{i:n}(t) = \frac{n!}{(i-1)!\,(n-i)!} \left\{ f^{KwG}(t;a,b\theta) \sum_{j=0}^{\infty} C_j \left\{F^{KwG}(t;a,b)\right\}^j \right\} \sum_{l=0}^{i-1} (-1)^l \binom{i-1}{l}$$

$$[\bar{F}^{KwG}(t;a,b)]^{\theta(n+l-i)} \sum_{k=0}^{\infty} d_{n+l-i,k} [F^{KwG}(t;a,b)]^k$$

$$= \frac{n!}{(i-1)!\,(n-i)!} \sum_{l=0}^{i-1} (-1)^l \binom{i-1}{l} f^{KwG}(t;a,b\theta) [\bar{F}^{KwG}(t;a,b)]^{\theta(n+l-i)}$$

$$\sum_{j,k=0}^{\infty} C_j\, d_{n+l-i,k} \left\{F^{KwG}(t;a,b)\right\}^{j+k}$$

$$= f^{KwG}(t;a,b\theta) [\bar{F}^{KwG}(t;a,b)]^{\theta(n+l-i)} \sum_{j,k=0}^{\infty} \phi_{j,k} \left\{F^{KwG}(t;a,b)\right\}^{j+k}$$

$$= f^{KwG}(t;a,b) [\bar{F}^{KwG}(t;a,b)]^{\theta-1} [\bar{F}^{KwG}(t;a,b)]^{\theta(n+l-i)} \sum_{j,k=0}^{\infty} \phi_{j,k} \left\{F^{KwG}(t;a,b)\right\}^{j+k}$$

$$= f^{KwG}(t;a,b)\, [\bar{F}^{KwG}(t;a,b)]^{\theta(n+l-i+1)-1} \sum_{j,k=0}^{\infty} \phi_{j,k} \left\{F^{KwG}(t;a,b)\right\}^{j+k} \qquad (18)$$

Where $$\phi_{j,k} = n\, C_j\, d_{n+l-i,k} \binom{n-1}{i-1} \sum_{l=0}^{i-1} (-1)^l \binom{i-1}{l}$$

**Remark1**. Equations (10)-(18) reveal that the density functions of the $GMOKw-G$ distribution and that of its order statistics can be expressed as different $Kw-G$ (Cordeiro and de Castro, 2011) density multiplied by an infinite power series of cdf of $Kw-G$ (Cordeiro and de Castro, 2011) and also as negative binomial mixture of $Kw-G$ (Cordeiro and de Castro, 2011) distributions under Lehman alternatives. These results are expected to play important roles to obtain explicit expressions for the moments and moment generating function (mgf) of the $GMOKw-G$ distribution and those of its order statistics in a general set up and for special models using the corresponding results of $Kw-G$ (Cordeiro and de Castro, 2011) distributions.

**5.3 Probability weighted moments**

The probability weighted moments (PWMs), first proposed by Greenwood et al. (1979), are expectations of certain functions of a random variable whose mean exists. The $(p,q,r)^{th}$ PWM of $T$ is defined by

$$\Gamma_{p,q,r} = \int_{-\infty}^{\infty} t^{p} [F(t)]^{q} [1-F(t)]^{r} f(t)\,dt$$

From equations (10), (14) and (15), the $s^{th}$ moment of $T$ for $\alpha \in (0,1)$ and $\alpha > 1$, can be written either as

$$E(T^{s}) = \int_{-\infty}^{+\infty} t^{s}\, a b\, g(t)\, G(t)^{a-1}[1-G(t)^{a}]^{b-1} \sum_{j=0}^{\infty} A_{j}\, [\{1-G(t)^{a}\}^{b}]^{j+\theta-1}\,dt$$

$$= \sum_{j=0}^{\infty} A_{j} \int_{-\infty}^{+\infty} t^{s} [\{1-G(t)^{a}\}^{b}]^{j+\theta-1}\, a b\, g(t) G(t)^{a-1}[1-G(t)^{a}]^{b-1}\,dt$$

$$= \sum_{j=0}^{\infty} A_{j}\, \Gamma_{s,0,\,j+\theta-1}$$

or

$$E(T^{s}) = \int_{-\infty}^{+\infty} t^{s}\, a b\, g(t) G(t)^{a-1}[1-G(t)^{a}]^{b\theta-1} \sum_{j=0}^{\infty} \sum_{k=0}^{j} \mathrm{B}_{j,k}\, [\,1-[1-G(t)^{a}]^{b}\,]^{j-k}\,dt$$

$$= \sum_{j=0}^{\infty} \sum_{k=0}^{j} B_{j,k} \int_{-\infty}^{+\infty} t^{s} [1-[1-G(t)^{a}]^{b}]^{j-k}\, a b\, g(t) G(t)^{a-1}[1-G(t)^{a}]^{b\theta-1}\,dt$$

$$= \sum_{j=0}^{\infty} \sum_{k=0}^{j} B_{j,k} \int_{-\infty}^{+\infty} t^{s} [1-[1-G(t)^{a}]^{b}]^{j-k} [1-G(t)^{a}]^{b(\theta-1)}\, a b\, g(t) G(t)^{a-1}[1-G(t)^{a}]^{b-1}\,dt$$

$$= \sum_{j=0}^{\infty} \sum_{k=0}^{j} B_{j,k}\, \Gamma_{s,\,j-k,\,\theta-1}$$

or

$$E(T^{s}) = \int_{-\infty}^{+\infty} t^{s}\, a b\, g(t)\, G(t)^{a-1}[1-G(t)^{a}]^{b\theta-1} \sum_{j=0}^{\infty} C_{j}\, [1-[1-G(t)^{a}]^{b}]^{j}\,dt$$

$$= \sum_{j=0}^{\infty} C_{j} \int_{-\infty}^{+\infty} t^{s} [1-[1-G(t)^{a}]^{b}]^{j}\, a b\, g(t) G(t)^{a-1}[1-G(t)^{a}]^{b\theta-1}\,dt$$

$$= \sum_{j=0}^{\infty} C_{j} \int_{-\infty}^{+\infty} t^{s} [1-[1-G(t)^{a}]^{b}]^{j} [1-G(t)^{a}]^{b(\theta-1)}\, a b\, g(t) G(t)^{a-1}[1-G(t)^{a}]^{b-1}\,dt$$

$$= \sum_{j=0}^{\infty} C_{j}\, \Gamma_{s,\,j,\,\theta-1}$$

Where $\Gamma_{p,q,r} = \int_{-\infty}^{\infty} t^{p}\, \{1-[1-G(t)^{a}]^{b}\}^{q} \{[1-G(t)^{a}]^{b}\}^{r} [a b\, g(t)\, G(t)^{a-1}[1-G(t)^{a}]^{b-1}]\,dt$

is the PWM of $\mathrm{Kw\text{-}G}(a,b)$ distribution.

Therefore the moments of the $\mathrm{GMOKw\text{-}G}(\alpha,\theta,a,b)$ can be expresses in terms of the PWMs of $\mathrm{Kw\text{-}G}(a,b)$ (Cordeiro and de Castro, 2011). The PWM method can generally be used for estimating parameters quantiles of generalized distributions. These moments have low variance and no severe biases, and they compare favourably with estimators obtained by maximum likelihood.

Proceeding as above we can derive $s^{th}$ moment of the $i^{th}$ order statistic $T_{i:n}$, in a random sample of size *n* from $GMOKw-G$ for $\alpha \in (0,1)$ and $\alpha > 1$, on using equations (16) and (18) as

$$E(T^{s}{}_{i,n}) = \sum_{j,p=0}^{\infty} M_{j,p}\, \Gamma_{s,0,\ j+p+\theta(n+l-i+1)-1} \text{ and } E(T^{s}{}_{i,n}) = \sum_{j,k=0}^{\infty} \phi_{j,k}\, \Gamma_{s,\ j+k,\ \theta(n+l-i+1)-1}$$

Where $A_j$, $B_{j,k}$, $C_j$ and $M_{j,p}, \phi_{j,k}$ defined in section 5.1 and 5.2 respectively.

### 5.4 Moment generating function

The moment generating function of $GMOKw-G$ family can be easily expressed in terms of those of the exponentiated $Kw-G$ (Cordeiro and de Castro, 2011) distribution using the results of section 5.1. For example using equation (11) it can be seen that

$$M_T(s) = E[e^{sT}] = \int_{-\infty}^{\infty} e^{st} f(t)\,dt = \int_{-\infty}^{\infty} e^{st} \sum_{j=0}^{\infty} A'_j \frac{d}{dt}[\overline{F}^{KwG}(t;a,b)]^{j+\theta}\,dt$$

$$= \sum_{j=0}^{\infty} A'_j \int_{-\infty}^{\infty} e^{st} \frac{d}{dt}\left\{\overline{F}^{KwG}(t;a,b)\right\}^{j+\theta} dt = \sum_{j=0}^{\infty} A'_j M_X(s)$$

Where $M_X(s)$ is the mgf of a $Kw-G$ (Cordeiro and de Castro, 2011) distribution.

### 5.5 Rényi Entropy

The entropy of a random variable is a measure of uncertainty variation and has been used in various situations in science and engineering. The Rényi entropy is defined by

$$I_R(\delta) = (1-\delta)^{-1} \log\left(\int_{-\infty}^{\infty} f(t)^{\delta}\,dt\right)$$

where $\delta > 0$ and $\delta \neq 1$ For furthers details, see Song (2001). For $\alpha \in (0,1)$ using expansion (9), in (6) we can write

$$f^{GMOKwG}(t)^{\delta} = \left[\frac{\theta\alpha^{\theta}\, a b\, g(t) G(t)^{a-1}[1-G(t)^{a}]^{b\theta-1}}{[1-\overline{\alpha}[1-G(t)^{a}]^{b}]^{\theta+1}}\right]^{\delta}$$

$$= [\theta\alpha^{\theta}\, a b\, g(t) G(t)^{a-1}\{1-G(t)^{a}\}^{b\theta-1}]^{\delta}\ [1-\overline{\alpha}[1-G(t)^{a}]^{b}]^{-\delta(\theta+1)}$$

$$= \frac{\theta^{\delta}\alpha^{\delta\theta}[a\,b\,g(t)G(t)^{a-1}\{1-G(t)^{a}\}^{b\theta-1}]^{\delta}}{\Gamma[\delta(\theta+1)]}\sum_{j=0}^{\infty}(1-\alpha)^{j}\,\Gamma[\delta(\theta+1)+j]\,\frac{[\{1-G(t)^{a}\}^{b}]^{j}}{j!}$$

Thus for $\alpha \in (0,1)$, the Rényi entropy of GMOKw - G $(\alpha,\theta,a,b)$ can be obtained as

$$I_R(\delta) = (1-\delta)^{-1}\log\left(\sum_{j=0}^{\infty} L_j \int_{-\infty}^{\infty}[a\,b\,g(t)G(t)^{a-1}\{1-G(t)^{a}\}^{b\theta-1}]^{\delta}[\{1-G(t)^{a}\}^{b}]^{j}\,dt\right)$$

Where $L_j = L_j(\alpha) = \dfrac{\theta^{\delta}\alpha^{\delta\theta}(1-\alpha)^{j}\,\Gamma[\delta(\theta+1)+j]}{\Gamma[\delta(\theta+1)]\;\; j!}$

Again the density function (6) can be expressed as

$$f^{GMOKwG}(t)^{\delta} = \left\{\theta\,a\,b\,g(t)G(t)^{a-1}[1-G(t)^{a}]^{b\theta-1}\Big/\alpha\left[1-\frac{(\alpha-1)[1-[1-G(t)^{a}]^{b}]}{\alpha}\right]^{\theta+1}\right\}^{\delta}$$

For $\alpha > 1$, using expansion (9) we get,

$$= \left\{\frac{\theta\,a\,b\,g(t)G(t)^{a-1}[1-G(t)^{a}]^{b\theta-1}}{\alpha}\left[1-\frac{(\alpha-1)[1-[1-G(t)^{a}]^{b}]}{\alpha}\right]^{-(\theta+1)}\right\}^{\delta}$$

$$= \frac{[\theta\,a\,b\,g(t)G(t)^{a-1}[1-G(t)^{a}]^{b\theta-1}]^{\delta}}{\alpha^{\delta+j}\,\Gamma[\delta(\theta+1)]}\sum_{j=0}^{\infty}\Gamma[\delta(\theta+1)+j](\alpha-1)^{j}\,\frac{[1-[1-G(t)^{a}]^{b}]^{j}}{j!}$$

Thus for $\alpha > 1$, the Rényi entropy of $GMOKw-G$ also can be derived as

$$I_R(\delta) = (1-\delta)^{-1}\log\left(\sum_{j=0}^{\infty} R_j \int_{-\infty}^{\infty}[a\,b\,g(t)G(t)^{a-1}[1-G(t)^{a}]^{b\theta-1}]^{\delta}[1-[1-G(t)^{a}]^{b}]^{j}\,dt\right)$$ Where

$$R_j = R_j(\alpha) = \frac{\theta^{\delta}(\alpha-1)^{j}\,\Gamma[\delta(\theta+1)+j]}{\alpha^{\delta+j}\,\Gamma[\delta(\theta+1)]\;j!}$$

**5.6 Quantile function and random sample generation**

We shall now present a formula for generating $GMOKw-G$ random variable by using inversion method by inverting the cdf or the survival function.

$$\overline{F}^{GMOKwG}(t) = \left[\frac{\alpha[1-G(t)^{a}]^{b}}{1-\overline{\alpha}[1-G(t)^{a}]^{b}}\right]^{\theta} \Rightarrow \frac{\alpha[1-G(t)^{a}]^{b}}{1-\overline{\alpha}[1-G(t)^{a}]^{b}} = \overline{F}^{GMOKwG}(t)^{1/\theta} \Rightarrow \alpha[1-G(t)^{a}]^{b}$$

$$= \overline{F}^{GMOKwG}(t)^{1/\theta}[1-\overline{\alpha}\{1-G(t)^{a}\}^{b}] \Rightarrow [1-G(t)^{a}]^{b}\{\alpha+\overline{\alpha}\,\overline{F}^{GMOKwG}(t)^{1/\theta}\} = \overline{F}^{GMOKwG}(t)^{1/\theta}$$

$$\Rightarrow [1-G(t)^{a}]^{b} = \frac{\overline{F}^{GMOKwG}(t)^{1/\theta}}{\alpha+\overline{\alpha}\,\overline{F}^{GMOKwG}(t)^{1/\theta}} \Rightarrow 1-G(t)^{a} = \left(\frac{\overline{F}^{GMOKwG}(t)^{1/\theta}}{\alpha+\overline{\alpha}\,\overline{F}^{GMOKwG}(t)^{1/\theta}}\right)^{1/b}$$

$$\Rightarrow G(t)^a = 1-\left(\frac{\bar{F}^{GMOKwG}(t)^{1/\theta}}{\alpha+\bar{\alpha}\,\bar{F}^{GMOKwG}(t)^{1/\theta}}\right)^{1/b} \Rightarrow G(t) = \left(1-\left(\frac{\bar{F}^{GMOKwG}(t)^{1/\theta}}{\alpha+\bar{\alpha}\,\bar{F}^{GMOKwG}(t)^{1/\theta}}\right)^{1/b}\right)^{1/a}$$

$$\Rightarrow t = G^{-1}\left(1-\left(\frac{[1-F^{GMOKwG}(t)]^{1/\theta}}{\alpha+\bar{\alpha}[1-F^{GMOKwG}(t)]^{1/\theta}}\right)^{1/b}\right)^{1/a}$$

to generate an random variable from $GMOKw-G$ first generate a $u \sim U(0,1)$ then use

$$\Rightarrow t = G^{-1}\left(1-\left(\frac{[1-u]^{1/\theta}}{\alpha+\bar{\alpha}[1-u]^{1/\theta}}\right)^{1/b}\right)^{1/a} \qquad (19)$$

The $p^{th}$ Quantile $t_p$ for $GMOKw-G$ can be easily obtained from (19) as

$$\Rightarrow t_p = G^{-1}\left[1-\left(\frac{[1-p]^{1/\theta}}{\alpha+\bar{\alpha}[1-p]^{1/\theta}}\right)^{1/b}\right]^{1/a}$$

For example, let the base line distribution be exponential with parameter $\lambda > 0$, having pdf and cdf as $g(t:\lambda) = \lambda e^{-\lambda t}, t > 0$ and $G(t:\lambda) = 1-e^{-\lambda t}$, respectively. Therefore the $p^{th}$ Quantile $t_p$ of $GMOKw-E$ is given by

$$t_p = -\frac{1}{\lambda}\log\left[1-\left[1-\left(\frac{[1-p]^{1/\theta}}{\alpha+\bar{\alpha}[1-p]^{1/\theta}}\right)^{1/b}\right]^{1/a}\right]$$

### 5.7 Asymptotes and shapes

Here we investigate the asymptotic shapes of the proposed family following the methods followed in Alizadeh *et al*., (2015).

**Proposition 2.** The asymptotes of equations (6), (7) and (8) as $t \to 0$ are given by

$$f(t) \sim \frac{\theta\, a\, b\, g(t) G(t)^{a-1}}{\alpha} \qquad \text{as } G(t) \to 0$$

$$F(t) \sim 0 \qquad \text{as } G(t) \to 0$$

$$h(t) \sim \frac{\theta\, a\, b\, g(t) G(t)^{a-1}}{\alpha} \qquad \text{as } G(t) \to 0$$

**Proposition 3.** The asymptotes of equations (6), (7) and (8) as $t \to \infty$ are given by

$$f(t) \sim \theta \alpha^{\theta} a\, b\, g(t)[1-G(t)^a]^{b\theta-1} \qquad \text{as } t \to \infty$$

$$F(t) \sim 1-\alpha^{\theta}[1-G(t)^{a}]^{b\theta} \qquad \text{as} \quad t \to \infty$$

$$h(t) \sim \theta\, a\, b\, g(t)[1-G(t)^{a}]^{-1} \qquad \text{as} \quad t \to \infty$$

The shapes of the density and hazard rate functions can be described analytically. The critical points of the $GMOKw-G$ density function are the roots of the equation:

$$\frac{d\log[f(t)]}{dt} = \frac{g'(t)}{g(t)} + (a-1)\frac{g(t)}{G(t)} + a(1-b\theta)\frac{g(t)G(t)^{a-1}}{1-G(t)^{a}}$$

$$-(\theta+1)\,\frac{\bar{\alpha}\, a\, b\, g(t)\, G(t)^{a-1}[1-G(t)^{a}]^{b-1}}{1-\bar{\alpha}[1-G(t)^{a}]^{b}} \qquad (20)$$

There may be more than one root to (20). If $t = t_0$ is a root of (20) then it corresponds to a local maximum, a local minimum or a point of inflexion depending on whether $\lambda(t_0) < 0,\ \lambda(t_0) > 0\ \ or\ \ \lambda(t_0) = 0$ where $\lambda(t) = \frac{d^2}{dt^2}\log[f(t)]$

$$\lambda(t) = \frac{g(t)g''(t)-[g'(t)]^2}{g(t)^2} + (a-1)\,\frac{G(t)g'(t)-g(t)^2}{G(t)}$$

$$+a(1-b\theta)[\,\frac{g'(t)G(t)^{a-1}}{1-G(t)^{a}} + \frac{(a-1)\,g(t)^2\, G(t)^{a-2}}{1-G(t)^{a}} + \frac{a\, g(t)^2\, G(t)^{2a-2}}{[1-G(t)^{a}]^2}]$$

$$-\frac{(\theta+1)\bar{\alpha}\, ab[g'(t)G(t)^{a-1}[1-G(t)^{a}]^{b-1}]}{1-\bar{\alpha}[1-G(t)^{a}]^{b}} - \frac{(\theta+1)\bar{\alpha}\, ab[(a-1)\,g(t)^2\, G(t)^{a-2}[1-G(t)^{a}]^{b-1}]}{1-\bar{\alpha}[1-G(t)^{a}]^{b}}$$

$$+\frac{(\theta+1)\bar{\alpha}\, ab[a(b-1)\,g(t)^2\, G(t)^{2a-2}[1-G(t)^{a}]^{b-2}]}{1-\bar{\alpha}[1-G(t)^{a}]^{b}} + (\theta+1)\left[\frac{\bar{\alpha}\, ab\, g(t)G(t)^{a-1}[1-G(t)^{a}]^{b-1}}{1-\bar{\alpha}[1-G(t)^{a}]^{b}}\right]^2$$

$$= \frac{g(t)g''(t)-[g'(t)]^2}{g(t)^2} + (a-1)\,\frac{G(t)g'(t)-g(t)^2}{G(t)}$$

$$+a(1-b\theta)[\,\frac{g'(t)G(t)^{a-1}}{1-G(t)^{a}} + \frac{(a-1)\,g(t)^2\, G(t)^{a-2}}{1-G(t)^{a}} + \frac{a\, g(t)^2\, G(t)^{2a-2}}{[1-G(t)^{a}]^2}]$$

$$-\frac{(\theta+1)\bar{\alpha}\, g'(t)}{\theta\, g(t)}[1-G(t)^{a}]^{b}\, h^{GMOKwG}(t) - \frac{(\theta+1)}{\theta}\,\bar{\alpha}\,(a-1)\,g(t)G(t)^{-1}[1-G(t)^{a}]^{b}\, h^{GMOKwG}(t)$$

$$+\frac{(\theta+1)}{\theta}\bar{\alpha}\, a(b-1)\, g(t)G(t)^{a-1}(t)\,\{1-G(t)^{a}\}^{b-1}\, h^{GMOKwG}(t) \; + \frac{(\theta+1)}{\theta^2}\{[1-G(t)^{a}]^{b}\bar{\alpha}\, h^{GMOKwG}(t)\}^2$$

The critical points of $h(t)$ are the roots of the equation

$$\frac{d\log[h(t)]}{dt} = \frac{g'(t)}{g(t)} + (a-1)\frac{g(t)}{G(t)} + \frac{a\, g(t)G(t)^{a-1}}{1-G(t)^{a}} - \frac{\bar{\alpha}\, a\, b\, g(t)G(t)^{a-1}[1-G(t)^{a}]^{b-1}}{1-\bar{\alpha}[1-G(t)^{a}]^{b}} \qquad (21)$$

There may be more than one root to (21). If $t = t_0$ is a root of (21) then it corresponds to a local maximum, a local minimum or a point of inflexion depending on whether

$\gamma(t_0)<0,\ \gamma(t_0)>0\ or\ \gamma(t_0)=0$, where $\gamma(t)=\frac{d^2}{dt^2}\log[h(t)]$

$$\gamma(t)=\frac{g(t)g''(t)-[g'(t)]^2}{g(t)^2}+(a-1)\frac{G(t)g'(t)-g(t)^2}{G(t)}$$

$$+a\left[\frac{g'(t)G(t)^{a-1}}{1-G(t)^a}+\frac{(a-1)g(t)^2G(t)^{a-2}}{1-G(t)^a}+\frac{a\,g(t)^2G(t)^{2a-2}}{[1-G(t)^a]^2}\right]$$

$$-\frac{\bar{\alpha}ab[g'(t)G(t)^{a-1}[1-G(t)^a]^{b-1}]}{1-\bar{\alpha}[1-G(t)^a]^b}-\frac{\bar{\alpha}ab[(a-1)g(t)^2G(t)^{a-2}[1-G(t)^a]^{b-1}]}{1-\bar{\alpha}[1-G(t)^a]^b}$$

$$+\frac{\bar{\alpha}ab[a(b-1)g(t)^2G(t)^{2a-2}[1-G(t)^a]^{b-2}]}{1-\bar{\alpha}[1-G(t)^a]^b}+\left[\frac{\bar{\alpha}ab\,g(t)G(t)^{a-1}[1-G(t)^a]^{b-1}}{1-\bar{\alpha}[1-G(t)^a]^b}\right]^2$$

$$=\frac{g(t)g''(t)-[g'(t)]^2}{g(t)^2}+(a-1)\frac{G(t)g'(t)-g(t)^2}{G(t)}$$

$$+a\left[\frac{g'(t)G(t)^{a-1}}{1-G(t)^a}+\frac{(a-1)g(t)^2G(t)^{a-2}}{1-G(t)^a}+\frac{a\,g(t)^2G(t)^{2a-2}}{[1-G(t)^a]^2}\right]-\frac{\bar{\alpha}g'(t)}{\theta g(t)}[1-G(t)^a]^b h^{GMOKwG}(t)$$

$$-\frac{1}{\theta}\bar{\alpha}(a-1)g(t)G(t)^{-1}(t)\,[1-G(t)^a]^b h^{GMOKwG}(t)+\frac{1}{\theta}\bar{\alpha}\,a(b-1)g(t)G(t)^{a-1}(t)\,\{1-G(t)^a\}^{b-1}h^{GMOKwG}(t)$$

$$+\frac{1}{\theta^2}\{[1-G(t)^a]^b\bar{\alpha}\,h^{GMOKwG}(t)\}^2$$

### 5.8 Stochastic orderings

In this section we study the reliability properties and stochastic ordering of the $GMOKw-G$ distributions Stochastic ordering properties have applications in diverse fields such as economics, reliability, survival analysis, insurance, finance, actuarial and management sciences (Shaked and Shanthikumar, 2007).

Let *X* and *Y* be two random variables with cfds *F* and *G*, respectively, survival functions $\bar{F}=1-F$ and $\bar{G}=1-G$, and corresponding pdf's *f, g*. Then *X* is said to be smaller than *Y* in the likelihood ratio order ($X\le_{lr}Y$) if $f(t)/g(t)$ is decreasing in $t\ge 0$; stochastic order ($X\le_{st}Y$) if $\bar{F}(t)\le\bar{G}(t)$ for all $t\ge 0$; hazard rate order ($X\le_{hr}Y$) if $\bar{F}(t)/\bar{G}(t)$ is decreasing in $t\ge 0$; reversed hazard rate order ($X\le_{rhr}Y$) if $F(t)/G(t)$ is decreasing in $t\ge 0$. These four stochastic orders are related to each other, as

$$X\le_{rhr}Y\Leftarrow X\le_{lr}Y\Rightarrow X\le_{hr}Y\Rightarrow X\le_{st}Y \qquad (22)$$

**Theorem1**: Let $\mathrm{X}\sim\mathrm{GMOKw\text{-}G}(\alpha_1,\theta,a,b)$ and $\mathrm{Y}\sim\mathrm{GMOKw\text{-}G}(\alpha_2,\theta,a,b)$. If $\alpha_1<\alpha_2$, then $X\le_{lr}Y$ ($X\le_{hr}Y,\ X\le_{rhr}Y,\ X\le_{st}Y$).

Proof: $\dfrac{f(t)}{g(t)} = \dfrac{\theta \alpha_1^{\theta} ab\, g(t) G(t)^{a-1}[1-G(t)^a]^{b\theta-1} \big/ [1-\bar{\alpha}_1[1-G(t)^a]^b]^{\theta+1}}{\theta \alpha_2^{\theta} ab\, g(t) G(t)^{a-1}[1-G(t)^a]^{b\theta-1} \big/ [1-\bar{\alpha}_2[1-G(t)^a]^b]^{\theta+1}}$

$$= \left(\frac{\alpha_1}{\alpha_2}\right)^{\theta} \left[\frac{1-\bar{\alpha}_2[1-G(t)^a]^b}{1-\bar{\alpha}_1[1-G(t)^a]^b}\right]^{\theta+1}$$

Since $\alpha_1 < \alpha_2$

$$\frac{d}{dt}\left[\frac{f(t)}{g(t)}\right] = (\theta+1)\left(\frac{\alpha_1}{\alpha_2}\right)^{\theta}\left[\frac{1-\bar{\alpha}_2[1-G(t)^a]^b}{1-\bar{\alpha}_1[1-G(t)^a]^b}\right]^{\theta}$$

$$\times \frac{\{1-\bar{\alpha}_1[1-G(t)^a]^b\}[\bar{\alpha}_2 ab[1-G(t)^a]^{b-1} G(t)^{a-1} g(t)] - \{1-\bar{\alpha}_2[1-G(t)^a]^b\}}{\{1-\bar{\alpha}_1[1-G(t)^a]^b\}^2}$$

$$[[\bar{\alpha}_1 ab[1-G(t)^a]^{b-1} G(t)^{a-1} g(t)]]$$

$$= (\theta+1)\left(\frac{\alpha_1}{\alpha_2}\right)^{\theta}\left[\frac{1-\bar{\alpha}_2[1-G(t)^a]^b}{1-\bar{\alpha}_1[1-G(t)^a]^b}\right]^{\theta}$$

$$\times \frac{ab[1-G(t)^a]^{b-1} G(t)^{a-1} g(t)[\{1-\bar{\alpha}_1[1-G(t)^a]^b\}\bar{\alpha}_2 - \{1-\bar{\alpha}_2[1-G(t)^a]^b\}\bar{\alpha}_1]}{\{1-\bar{\alpha}_1[1-G(t)^a]^b\}^2}$$

$$= (\theta+1)\left(\frac{\alpha_1}{\alpha_2}\right)^{\theta} \frac{\{1-\bar{\alpha}_2[1-G(t)^a]^b\}^{\theta} ab[1-G(t)^a]^{b-1} G(t)^{a-1} g(t)(\bar{\alpha}_2 - \bar{\alpha}_1)}{\{1-\bar{\alpha}_1[1-G(t)^a]^b\}^{\theta+2}}$$

$$= (\theta+1)\left(\frac{\alpha_1}{\alpha_2}\right)^{\theta} (\alpha_1 - \alpha_2) \frac{\{1-\bar{\alpha}_2[1-G(t)^a]^b\}^{\theta} ab[1-G(t)^a]^{b-1} G(t)^{a-1} g(t)}{\{1-\bar{\alpha}_1[1-G(t)^a]^b\}^{\theta+2}}$$

Which is always less than 0.

Hence, $f(t)/g(t)$ is decreasing in *t*. That is $X \leq_{lr} Y$ .The remaining statements follow from the implications (22).

## 6. Estimation

### 6.1 Maximum likelihood method

The model parameters of the $GMOKw-G$ distribution can be estimated by maximum likelihood. Let $t=(t_1,t_2,...t_n)^T$ be a random sample of size $n$ from $GMOKw-G$ with parameter vector $\boldsymbol{\theta}=(\theta,\alpha,a,b,\boldsymbol{\beta}^T)^T$, where $\boldsymbol{\beta}=(\beta_1,\beta_2,...\beta_q)^T$ corresponds to the parameter vector of the baseline distribution *G*. Then the log-likelihood function for $\boldsymbol{\theta}$ is given by

$$\ell = \ell(\boldsymbol{\theta}) = n\log\theta + n\theta\log\alpha + n\log(ab) + \sum_{i=0}^{n}\log[g(t_i,\boldsymbol{\beta})] + (a-1)\sum_{i=0}^{n}\log[G(t_i,\boldsymbol{\beta})]$$

$$+(b\theta-1)\sum_{i=0}^{n}\log[1-G(t_i,\boldsymbol{\beta})^a] - (\theta+1)\sum_{i=1}^{n}\log[1-\bar{\alpha}[1-G(t_i,\boldsymbol{\beta})^a]^b] \qquad (23)$$

This log-likelihood function can not be solved analytically because of its complex form but it can be maximized numerically by employing global optimization methods available with software's like R, SAS, Mathematica or by solving the nonlinear likelihood equations obtained by differentiating (23).

By taking the partial derivatives of the log-likelihood function with respect to $\theta,\alpha,a,b$ and $\boldsymbol{\beta}$ we obtain the components of the score vector $U_{\boldsymbol{\theta}}=(U_\theta,U_\alpha,U_a,U_b,U_{\beta^T})^T$

$$U_\theta=\frac{\partial\ell}{\partial\theta}=\frac{n}{\theta}+n\log\alpha+b\sum_{i=0}^{n}\log[1-G(t_i,\boldsymbol{\beta})^a]-\sum_{i=1}^{n}\log[1-\bar{\alpha}[1-G(t_i,\boldsymbol{\beta})^a]^b]$$

$$U_\alpha=\frac{\partial\ell}{\partial\alpha}=\frac{n\theta}{\alpha}-(\theta+1)\sum_{i=0}^{n}\frac{[1-G(t_i,\boldsymbol{\beta})^a]^b}{1-\bar{\alpha}\,[1-G(t_i,\boldsymbol{\beta})^a]^b}$$

$$U_a=\frac{\partial\ell}{\partial a}=\frac{n}{a}+\sum_{i=0}^{n}\log[G(t_i,\boldsymbol{\beta})]+(1-b\theta)\sum_{i=0}^{n}\frac{G(t_i,\boldsymbol{\beta})^a\log[G(t_i,\boldsymbol{\beta})]}{1-G(t_i,\boldsymbol{\beta})^a}$$
$$-(\theta+1)\sum_{i=0}^{n}\frac{b\bar{\alpha}[1-G(t_i,\boldsymbol{\beta})^a]^{b-1}G(t_i,\boldsymbol{\beta})^a\log[G(t_i,\boldsymbol{\beta})]}{1-\bar{\alpha}\,[1-G(t_i,\boldsymbol{\beta})^a]^b}$$

$$U_b=\frac{\partial\ell}{\partial b}=\frac{n}{b}+\theta\sum_{i=0}^{n}\log[1-G(t_i,\boldsymbol{\beta})^a]+(\theta+1)\sum_{i=0}^{n}\frac{\bar{\alpha}[1-G(t_i,\boldsymbol{\beta})^a]^b\log[1-G(t_i,\boldsymbol{\beta})^a]}{1-\bar{\alpha}[1-G(t_i,\boldsymbol{\beta})^a]^b}$$

$$U_\beta=\frac{\partial\ell}{\partial\boldsymbol{\beta}}=\sum_{i=0}^{n}\frac{g^{(\boldsymbol{\beta})}(t_i,\boldsymbol{\beta})}{g(t_i,\boldsymbol{\beta})}+(a-1)\sum_{i=0}^{n}\frac{G^{(\boldsymbol{\beta})}(t_i,\boldsymbol{\beta})}{G(t_i,\boldsymbol{\beta})}+(1-b\theta)\sum_{i=0}^{n}\frac{aG(t_i,\boldsymbol{\beta})^{a-1}G^{(\boldsymbol{\beta})}(t_i,\boldsymbol{\beta})}{1-G(t_i,\boldsymbol{\beta})^a}$$
$$-(\theta+1)\sum_{i=0}^{n}\frac{b\bar{\alpha}[1-G(t_i,\boldsymbol{\beta})^a]^{b-1}aG(t_i,\boldsymbol{\beta})^{a-1}G^{(\boldsymbol{\beta})}(t_i,\boldsymbol{\beta})}{1-\bar{\alpha}\,[1-G(t_i,\boldsymbol{\beta})^a]^b}$$

Setting these equations to zero $U_{\boldsymbol{\theta}}=(U_\theta,U_\alpha,U_a,U_b,U_{\beta^T})^T=0$ and solving them simultaneously yields the maximum likelihood estimate (MLE) $\hat{\boldsymbol{\theta}}=(\hat{\theta},\hat{\alpha},\hat{a},\hat{b},\hat{\boldsymbol{\beta}}^T)^T$ of $\boldsymbol{\theta}=(\theta,\alpha,a,b,\boldsymbol{\beta}^T)^T$.

**6.2 Asymptotic standard error and confidence interval for the mles:**

The asymptotic variance-covariance matrix of the MLEs of parameters can obtained by inverting the Fisher information matrix $\mathrm{I}(\boldsymbol{\theta})$ which can be derived using the second partial derivatives of the log-likelihood function with respect to each parameter. The $i\,j^{th}$ elements of $\mathrm{I}_n(\boldsymbol{\theta})$ are given by

$$\mathrm{I}_{ij}=-E\left(\frac{\partial^2 l(\boldsymbol{\theta})}{\partial\theta_i\,\partial\theta_j}\right),\qquad i,\ j=1,2,\cdots,3+q$$

The exact evaluation of the above expectations may be cumbersome. In practice one can estimate $\mathrm{I}_n(\boldsymbol{\theta})$ by the observed Fisher's information matrix $\hat{\mathrm{I}}_n(\hat{\boldsymbol{\theta}})$ is defined as:

$$\hat{\mathrm{I}}_{ij} \approx \left( -\frac{\partial^2 l(\boldsymbol{\theta})}{\partial\theta_i \partial\theta_j} \right)_{\boldsymbol{\theta}=\hat{\boldsymbol{\theta}}} , \qquad i,\ j = 1,\ 2, \cdots, 3+q$$

Using the general theory of MLEs under some regularity conditions on the parameters as $n \to \infty$ the asymptotic distribution of $\sqrt{n}\,(\hat{\boldsymbol{\theta}} - \boldsymbol{\theta})$ is $N_k(0, V_n)$ where $V_n = (v_{jj}) = \mathrm{I}_n^{-1}(\boldsymbol{\theta})$. The asymptotic behaviour remains valid if $V_n$ is replaced by $\hat{V}_n = \hat{\mathrm{I}}^{-1}(\hat{\boldsymbol{\theta}})$. This result can be used to provide large sample standard errors and also construct confidence intervals for the model parameters. Thus an approximate standard error and $(1-\gamma/2)100\%$ confidence interval for the mle of $j^{\text{th}}$ parameter $\theta_j$ are respectively given by $\sqrt{\hat{v}_{jj}}$ and $\hat{\theta}_j \pm Z_{\gamma/2}\sqrt{\hat{v}_{jj}}$ , where $Z_{\gamma/2}$ is the $\gamma/2$ point of standard normal distribution.

As an illustration on the MLE method its large sample standard errors, confidence interval in the case of $GMOKw - E(\theta, \alpha, a, b, \lambda)$ is discussed in an appendix.

### 6.3 Estimation by method of moments

Here an alternative method to estimation of the parameters is discussed. Since the moments are not in closed form, the estimation by the usual method of moments may not be tractable. Therefore following (Barreto-Souzai *et al.*, 2013) we get

$$E\,[\ 1-\bar{\alpha}[1-G(t)^a]^b\ ]$$

$$= \int_{-\infty}^{\infty} \{[\ 1-\bar{\alpha}[1-G(t)^a]^b\ ]\} \frac{\alpha^\theta\,\theta\, a\, b\, g(t) G(t)^{a-1}[1-G(t)^a]^{b\theta-1}}{[1-\bar{\alpha}[1-G(t)^a]^b\,]^{\theta+1}}\, dt$$

$$= \int_{-\infty}^{\infty} \frac{\alpha^\theta\,\theta\, a\, b\, g(t) G(t)^{a-1}[1-G(t)^a]^{b\theta-1}}{[1-\bar{\alpha}[1-G(t)^a]^b\,]^{\theta}}\, dt$$

Let, $u = 1-\bar{\alpha}\,[1-G(t)^a]^b$ then $du = \bar{\alpha}\, a\, b\, g(t) G(t)^{a-1}[1-G(t)^a]^{b-1}\, dt$

$$= \frac{\theta\alpha^\theta}{\bar{\alpha}} \int_{\alpha}^{1} \frac{\left((1-u)/\bar{\alpha}\right)^{(\theta-1)}}{u^\theta}\, du = \theta\left(\frac{\alpha}{\bar{\alpha}}\right)^{\theta} \int_{\alpha}^{1} u^{-\theta}\,(1-u)^{\theta-1}\, du$$

$$= \theta\left(\frac{\alpha}{\bar{\alpha}}\right)^{\theta} \int_{\alpha}^{1} u^{-1}\,(1/u-1)^{\theta-1}\, du$$

Let $1/u - 1 = v$

$\Rightarrow -du = u^2 dv$

$$= \theta\left(\frac{\alpha}{\bar{\alpha}}\right)^{\theta} \int_{0}^{(1/\alpha)-1} (1+v)^{-1}\, v^{\theta-1}\, dv = \theta\left(\frac{\alpha}{\bar{\alpha}}\right)^{\theta} \sum_{i\ge 0} \int_{0}^{(1/\alpha)-1} (-1)^i v^{\theta+i-1} dv$$

$$=\theta\left(\frac{\alpha}{\bar{\alpha}}\right)^{\theta}\sum_{i\geq 0}(-1)^{i}\left[\frac{v^{\theta+i}}{\theta+i}\right]_{0}^{\bar{\alpha}/\alpha}$$

$$=\theta\left(-\frac{\alpha}{\bar{\alpha}}\right)^{\theta}\sum_{i\geq 0}\left[\frac{(-\bar{\alpha}/\alpha)^{\theta+i}}{\theta+i}\right]$$

$$=(-1)^{\theta}\,\theta\sum_{i\geq 0}\left[\frac{(-\bar{\alpha}/\alpha)^{i}}{\theta+i}\right],\ \alpha>1/2$$

For $\theta=1$ we get the corresponding result of $\mathrm{MOKw\text{-}G}(\alpha,a,b)$ (Handique and Chakraborty, 2015).

$$=\left(\frac{\alpha}{\bar{\alpha}}\right)\sum_{i\geq 0}\left[\frac{(-\bar{\alpha}/\alpha)^{i+1}}{i+1}\right]=\left(\frac{\alpha}{\bar{\alpha}}\right)\left[\log\left(1+\frac{\bar{\alpha}}{\alpha}\right)\right]$$

$$=\left(\frac{\alpha}{\bar{\alpha}}\right)\left[-\log\left(\frac{1}{\alpha}\right)\right]=-\left(\frac{\alpha}{1-\alpha}\right)\log(\alpha)$$

$$E\{[\ 1-\bar{\alpha}[1-G(t)^{a}]^{b}\ ]\}^{\upsilon}=\int_{-\infty}^{\infty}\{[\ 1-\bar{\alpha}[1-G(t)^{a}]^{b}\ ]\}^{\upsilon}\frac{\alpha^{\theta}\,\theta\,a\,b\,g(t)G(t)^{a-1}[1-G(t)^{a}]^{b\theta-1}}{[1-\bar{\alpha}[1-G(t)^{a}]^{b}\ ]^{\theta+1}}\,dt$$

$$=\int_{-\infty}^{\infty}\frac{\alpha^{\theta}\,\theta\,a\,b\,g(t)G(t)^{a-1}[1-G(t)^{a}]^{b\theta-1}}{[1-\bar{\alpha}[1-G(t)^{a}]^{b}\ ]^{\theta+1-\upsilon}}\,dt$$

Let, $u=1-\bar{\alpha}\,[1-G(t)^{a}]^{b}$ then $du=\bar{\alpha}\,a\,b\,g(t)G(t)^{a-1}[1-G(t)^{a}]^{b-1}\,dt$

$$=\frac{\theta\,\alpha^{\theta}}{\bar{\alpha}}\int_{\alpha}^{1}\frac{\left((1-u)/\bar{\alpha}\right)^{(\theta-1)}}{u^{\theta+1-\upsilon}}\,du$$

$$=\theta\left(\frac{\alpha}{\bar{\alpha}}\right)^{\theta}\int_{\alpha}^{1}u^{\upsilon-\theta-1}(1-u)^{\theta-1}\,du$$

$$=\theta\left(\frac{\alpha}{\bar{\alpha}}\right)^{\theta}B_{\alpha}(\upsilon-\theta,\theta)\quad\text{for }\alpha<1,\upsilon>\theta$$

Where $B_{x}(m,n)=\int_{x}^{1}u^{m-1}(1-u)^{n-1}du$ is the incomplete beta function.

For $\theta=1$ we get the corresponding result of $\mathrm{MOKw\text{-}G}(\alpha,a,b)$

$$=\frac{\alpha(1-\alpha^{v-1})}{\bar{\alpha}(v-1)}.$$ (Handique and Chakraborty, 2015).

Therefore we have

$$E[\,[1-\overline{\alpha}[1-G(t)^a]^b]^v\,] = \begin{cases} (-1)^{\theta}\,\theta \sum_{i\geq 0}\left[\dfrac{(-\overline{\alpha}/\alpha)^i}{\theta+i}\right] & \text{for } \theta<1,\ \alpha>1/2, v=1 \\ \theta\left(\dfrac{\alpha}{\overline{\alpha}}\right)^{\theta} B_{\alpha}(\upsilon-\theta,\theta)\ , & \text{for } \alpha<1, \upsilon>\theta, v=2,3,... \end{cases} \tag{24}$$

For a random sample $t_1, t_2, ... t_n$ from a population with survival function (7), the model parameters can be estimated using (24) by solving the equations

$$\frac{1}{n}\sum_{i=1}^{n}[[1-\overline{\alpha}[1-G(t_i)^a]^b]^v] = \begin{cases} (-1)^{\theta}\,\theta \sum_{i\geq 0}\left[\dfrac{(-\overline{\alpha}/\alpha)^i}{\theta+i}\right] & \text{for } \theta<1,\ \alpha<1, v=1 \\ \theta\left(\dfrac{\alpha}{\overline{\alpha}}\right)^{\theta} B_{\alpha}(\upsilon-\theta,\theta)\ , & \text{for } \alpha<1, \upsilon>\theta, v=2,3,... \end{cases}$$

**6.4 Real life applications**

We consider one real life data to illustrate the suitability of the $GMOKw-W$ distribution presented in Section 4.2. The data set is a subset of data reported Bekker et al. (2000) which corresponds to the survival times (in years) of a group of patients given chemotherapy treatment alone. The data consisting of survival times (in years) for 46 patients are:

{0.047, 0.115, 0.121, 0.132, 0.164, 0.197, 0.203, 0.260, 0.282, 0.296, 0.334, 0.395, 0.458, 0.466, 0.501, 0.507, 0.529, 0.534, 0.540, 0.641, 0.644, 0.696, 0.841, 0.863, 1.099, 1.219, 1.271, 1.326, 1.447, 1.485, 1.553, 1.581, 1.589, 2.178, 2.343, 2.416, 2.444, 2.825, 2.830, 3.578, 3.658, 3.743, 3.978, 4.003, 4.033}

Here we have estimated the parameters by numerical maximization of loglikelihood function and provided their standard errors and 95% confidence intervals using large sample approach (see appendix). We have compared the $GMOKw-W$ distribution with it some of its sub models namely the Marshall-Olkin Kumaraswamy Weibull ($MOKw-W$) (Handique and Chakraborty, 2015), Kumaraswamy Weibull ($Kw-G$) and Marshall-Olkin Weibull ($MO-W$) and further with $BKw-W$ (Handique and Chakraborty, 2016) and $B-W$ models as well.

The best model is chosen as the one having lowest AIC (Akaike Information Criterion). It may be noted that $AIC = 2k - 2l$; where $k$ is the number of parameters in the statistical model and $l$ is the maximized value of the log-likelihood function under the considered model. The model with the lowest AIC value is considered as the best.

Likelihood Ratio Test for nested models:

The $\text{GMOKw-W}(a,b,\lambda,\alpha,\beta,\theta)$ distribution reduces to $\text{MOKw-W}(a,b,\lambda,\alpha,\beta)$ when $\theta=1$, to $\text{Kw-W}(a,b,\lambda,\beta)$ if $\alpha=\theta=1$, to $\text{MO-W}(\lambda,\beta,\theta)$ for $a=1,b=1,\theta=1$.

Here we have employed likelihood ratio criterion to test the following null hypothesis:

(i) $H_0:\theta=1$, that is the sample is from $\text{MOKw-W}(a,b,\lambda,\alpha,\beta)$

$H_1:\theta\neq 1$, that is the sample is $\text{GMOKw-W}(a,b,\lambda,\alpha,\beta,\theta)$ .

(ii) $H_0:\alpha=\theta=1$, that is the sample is from $\text{Kw-W}(a,b,\lambda,\beta)$

$H_1:\alpha\neq 1,\theta\neq 1$, that is the sample is $\text{GMOKw-W}(a,b,\lambda,\alpha,\beta,\theta)$ .

(iii) $H_0:a=1,b=1,\theta=1$, that is the sample is from $\text{MO-W}(\lambda,\alpha,\beta)$

$H_1:a\neq 1,b\neq 1,\theta\neq 1$, that is the sample is $\text{GMOKw-W}(a,b,\lambda,\alpha,\beta,\theta)$ .

Writing $\boldsymbol{\beta}=(a,b,\lambda,\alpha,\beta,\theta)$ the likelihood ratio test statistic is given by $\text{LR}=-2\ln\frac{L(\hat{\boldsymbol{\beta}}^{*};x)}{L(\hat{\boldsymbol{\beta}},x)}$, where $\hat{\boldsymbol{\beta}}^{*}$ is the restricted ML estimates under the null hypothesis $H_0$ and $\hat{\boldsymbol{\beta}}$ is the unrestricted ML estimates under the alternative hypothesis $H_1$. Under the null hypothesis $H_0$ the LR criterion follows Chi-square distribution with degrees of freedom (df) $(df_{alt}-df_{null})$. The null hypothesis is rejected for $p$-value less than 0.05.

In Table 1 the MLEs and standard errors and 95% confidence intervals (in parentheses) of the parameters for the fitted model along with the AIC values and Likelihood ratio statistic are presented. The $GMOKw-W$ model with the lowest values of the AIC is seen as the best model. Moreover LR test rejects all the three sub models in favour of the $GMOKw-W$.

The plots of the fitted densities and fitted cdf's along with the observed ones are displayed in Figures 3 and 4 indicate that the $GMOKw-W$ provides a best fit to the data considered here.

**Table 1:** MLEs, standard errors and 95% confidence intervals (in parentheses) along with log likelihood and AIC values for survival times of cancer patients data.

| **Distribution** | $MO-W$ | $Kw-W$ | $MOKw-W$ | $GMOKw-W$ | $B-W$ | $BKw-W$ |
|---|---|---|---|---|---|---|
| $\hat{a}$ | --- | 2.160<br>(1.584)<br>(-0.94, 5.26) | 0.737<br>(0.938)<br>(-1.10, 2.58) | 0.518<br>(0.011)<br>(0.49, 0.54) | --- | 6.525<br>(0.079)<br>(6.37, 6.67) |
| $\hat{b}$ | --- | 0.208<br>(0.204)<br>(-0.19, 0.61) | 0.322<br>(0.504)<br>(-0.67, 1.31) | 0.244<br>(0.113)<br>(0.02, 0.47) | --- | 0.317<br>(0.323)<br>(-0.32, 0.95) |
| $\hat{\lambda}$ | 1.229<br>(0.247)<br>(0.74, 1.71) | 4.521<br>(3.900)<br>(-3.12, 12.16) | 0.562<br>(2.102)<br>(-3.56, 4.68) | 0.111<br>(0.004)<br>(0.10, 0.12) | 0.467<br>(0.807)<br>(-1.11, 2.05) | 5.183<br>(0.004)<br>(5.17, 5.19) |
| $\hat{\beta}$ | 0.449<br>(0.451)<br>(-0.43, 1.33) | 0.836<br>(0.196)<br>(0.45, 1.22) | 1.598<br>(1.658)<br>(-1.65, 4.85) | 4.112<br>(0.005)<br>(4.10,4.12) | 0.695<br>(0.555)<br>(-0.39, 1.78) | 1.388<br>(0.002)<br>(1.38, 1.39) |
| $\hat{\alpha}$ | 0.418<br>(0.308)<br>(-0.19, 1.02) | --- | 0.240<br>(0.367)<br>(-0.47, 0.96) | 0.004<br>(0.002)<br>(0.00008, 0.008) | 2.024<br>(2.856)<br>(-3.57, 7.62) | 0.223<br>(0.054)<br>(0.12, 0.33) |
| $\hat{\theta}$ | --- | --- | --- | 0.239<br>(0.059)<br>(0.12, 0.4) | 3.418<br>(6.538)<br>(-9.39,16.2) | 0.224<br>(0.189)<br>(-0.15, 0.59) |
| $(l_{\max})$ | -57.87 | -57.72 | -57.81 | **-53.82** | -58.07 | -55.46 |
| AIC | 121.74 | 123.44 | 125.62 | **119.64** | 124.14 | 122.92 |
| LR | 8.1 | 7.8 | 7.98 | --- | --- | --- |
| ***p*-value** | **0.04399** | **0.02024** | **0.00473** | | | |

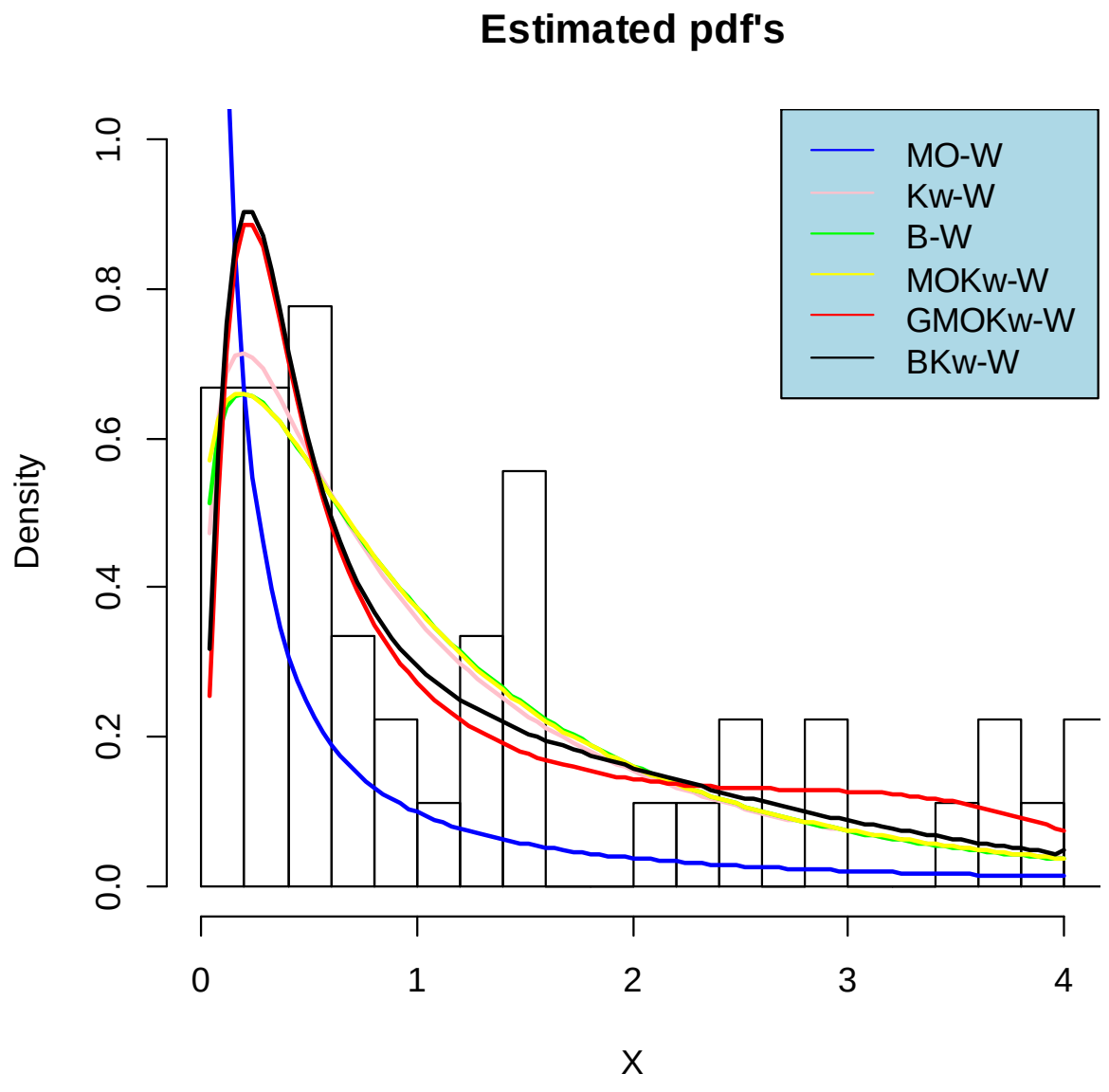

**Fig: 3** Plots of the observed histogram and estimated pdf's for the $MO-W$, $Kw-W$, $MOKw-W$ and $GMOKw-W$ models and $B-W$, $BKw-W$ models.

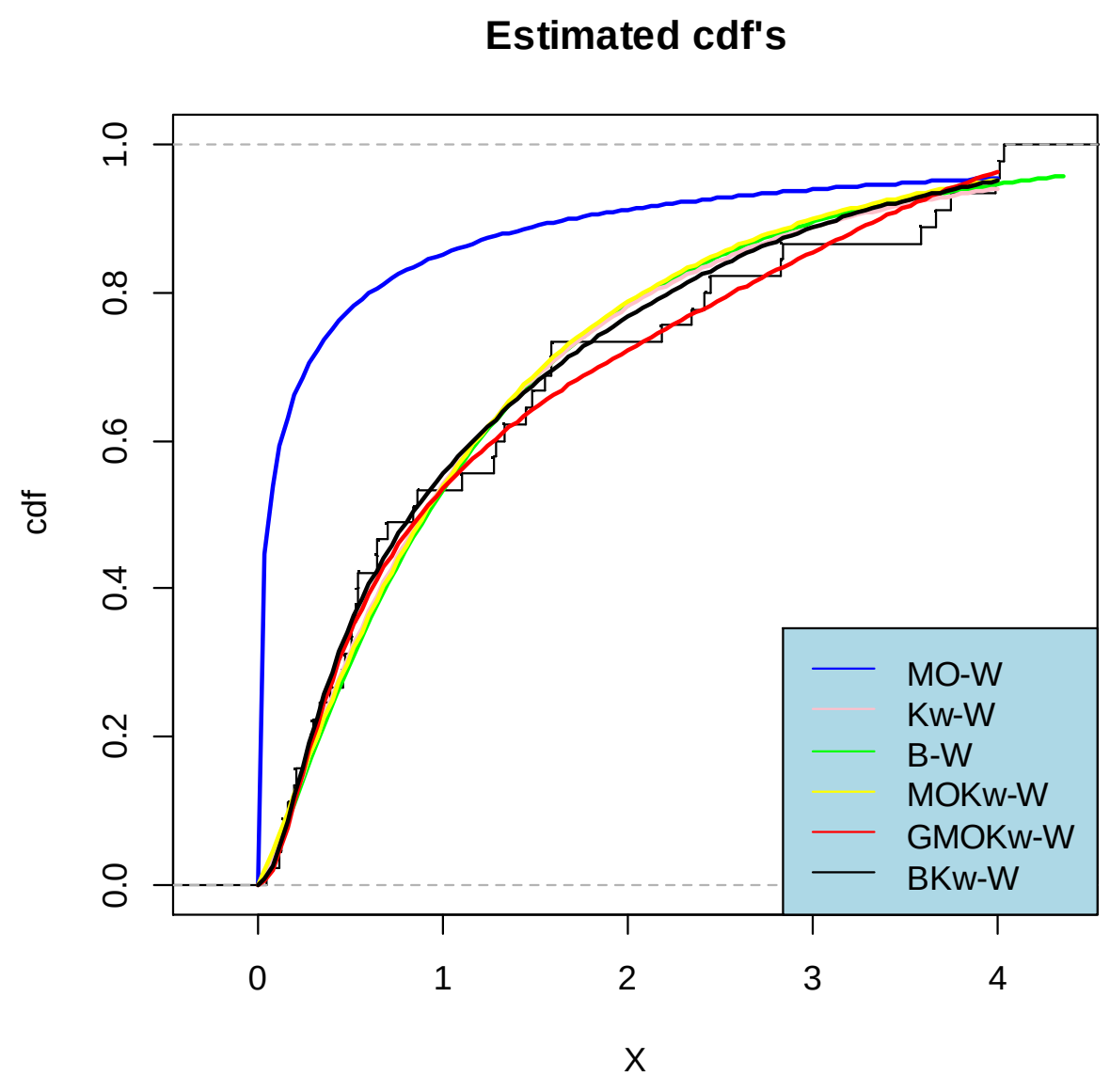

**Fig:4** Plots of the observed ogive and estimated cdf's for the $MO-W$, $Kw-W$, $MOKw-W$ and $GMOKw-W$ models and $B-W$, $BKw-W$ models .

### 7. Conclusion

Generalized Marshall-Olkin extended Kumaraswamy generalized family of distributions is introduced and some of its important properties and parameter estimation are studied. Applications with one real life data set has shown that the $GMOKw-W$ distribution is better than its sub models and also some other competitors with respect to the AIC values.

## Appendix: **Maximum likelihood estimation for** $GMOKw-E$

The pdf of the $GMOKw-E$ distribution is given by

$$f^{GMOKwE}(t) = \frac{\theta\alpha^{\theta}\, ab\lambda e^{-\lambda t}[1-e^{-\lambda t}]^{a-1}[1-[1-e^{-\lambda t}]^{a}]^{b\theta-1}}{[1-\bar{\alpha}[1-[1-e^{-\lambda t}]^{a}]^{b}]^{\theta+1}}$$

for $\theta > 0, \lambda > 0, \alpha > 0,\ a > 0, b > 0,\ t > 0$

For a random sample of size $n$ from this distribution, the log-likelihood function for the parameter vector $\boldsymbol{\theta} = (\theta, \alpha, a, b, \lambda)^T$ is given by

$$\ell = \ell(\boldsymbol{\theta}) = n\log\theta + n\theta\log\alpha + n\log(ab) + n\log\lambda - \lambda\sum_{i=0}^{n} t_i + (a-1)\sum_{i=0}^{n}\log(1-e^{-\lambda t_i})$$

$$+(b\theta-1)\sum_{i=0}^{n}\log[1-(1-e^{-\lambda t_i})^{a}] - (\theta+1)\sum_{i=0}^{n}\log[1-\bar{\alpha}[1-(1-e^{-\lambda t_i})^{a}]^{b}]$$

The components of the score vector $\boldsymbol{\theta} = (\theta, \alpha, a, b, \lambda)^T$ are

$$\frac{\partial \ell(\boldsymbol{\theta})}{\partial \theta} = \frac{n}{\theta} + n\log\alpha \; + b\sum_{i=0}^{n} \log[1-(1-e^{-\lambda t_i})^a] - \sum_{i=0}^{n} \log\,[1-\overline{\alpha}\,[1-(1-e^{-\lambda t_i})^a]^b]$$

$$\frac{\partial \ell(\boldsymbol{\theta})}{\partial \alpha} = \frac{n\theta}{\alpha} - (\theta+1)\sum_{i=0}^{n} \frac{[1-(1-e^{-\lambda t_i})^a]^b}{1-\overline{\alpha}\,[1-(1-e^{-\lambda t_i})^a]^b}$$

$$\frac{\partial \ell(\boldsymbol{\theta})}{\partial a} = \frac{n}{a} + \sum_{i=0}^{n} \log(1-e^{-\lambda t_i}) + (1-b\theta)\sum_{i=0}^{n} \frac{(1-e^{-\lambda t_i})^a \log(1-e^{-\lambda t_i})}{1-(1-e^{-\lambda t_i})^a}$$

$$-(\theta+1)\sum_{i=0}^{n} \frac{\overline{\alpha}\, b[1-(1-e^{-\lambda t_i})^a]^{b-1}(1-e^{-\lambda t_i})^a \log\,[(1-e^{-\lambda t_i})\,]}{1-\overline{\alpha}\,[1-(1-e^{-\lambda t_i})^a]^b}$$

$$\frac{\partial \ell(\boldsymbol{\theta})}{\partial b} = \frac{n}{b} + \theta\sum_{i=0}^{n} \log[1-(1-e^{-\lambda t_i})^a] + (\theta+1)\sum_{i=0}^{n} \frac{\overline{\alpha}\,[1-(1-e^{-\lambda t_i})^a]^b \log[1-(1-e^{-\lambda t_i})^a]}{1-\overline{\alpha}\,[1-(1-e^{-\lambda t_i})^a]^b}$$

$$\frac{\partial \ell(\boldsymbol{\theta})}{\partial \lambda} = \frac{n}{\lambda} - \sum_{i=0}^{n} t_i + (a-1)\sum_{i=0}^{n} \frac{\lambda e^{-\lambda t_i}}{1-e^{-\lambda t_i}} + (1-b\theta)\sum_{i=0}^{n} \frac{a\lambda(1-e^{-\lambda t_i})^{a-1} e^{-\lambda t_i}}{1-(1-e^{-\lambda t_i})^a}$$

$$-(\theta+1)\sum_{i=0}^{n} \frac{\overline{\alpha}\, ab\lambda\,[1-(1-e^{-\lambda t_i})^a]^{b-1}(1-e^{-\lambda t_i})^{a-1} e^{-\lambda t_i}}{1-\overline{\alpha}\,[1-(1-e^{-\lambda t_i})^a]^b}$$

The asymptotic variance covariance matrix for mles of the unknown parameters $\boldsymbol{\theta} = (\theta, \alpha, a, b, \lambda)$ of $GMOKw-E\ (\theta, \alpha, a, b, \lambda)$ distribution is estimated by

$$\hat{\mathrm{I}}_n{}^{-1}(\hat{\boldsymbol{\theta}}) = \begin{pmatrix} \mathrm{var}(\hat{\theta}) & \mathrm{cov}(\hat{\theta},\hat{\alpha}) & \mathrm{cov}(\hat{\theta},\hat{a}) & \mathrm{cov}(\hat{\theta},\hat{b}) & \mathrm{cov}(\hat{\theta},\hat{\lambda}) \\ \mathrm{cov}(\hat{\alpha},\hat{\theta}) & \mathrm{var}(\hat{\alpha}) & \mathrm{cov}(\hat{\alpha},\hat{a}) & \mathrm{cov}(\hat{\alpha},\hat{b}) & \mathrm{cov}(\hat{\alpha},\hat{\lambda}) \\ \mathrm{cov}(\hat{a},\hat{\theta}) & \mathrm{cov}(\hat{a},\hat{\alpha}) & \mathrm{var}(\hat{a}) & \mathrm{cov}(\hat{a},\hat{b}) & \mathrm{cov}(\hat{a},\hat{\lambda}) \\ \mathrm{cov}(\hat{b},\hat{\theta}) & \mathrm{cov}(\hat{b},\hat{\alpha}) & \mathrm{cov}(\hat{b},\hat{a}) & \mathrm{var}(\hat{b}) & \mathrm{cov}(\hat{b},\hat{\lambda}) \\ \mathrm{cov}(\hat{\lambda},\hat{\theta}) & \mathrm{cov}(\hat{\lambda},\hat{\alpha}) & \mathrm{cov}(\hat{\lambda},\hat{a}) & \mathrm{cov}(\hat{\lambda},\hat{b}) & \mathrm{var}(\hat{\lambda}) \end{pmatrix}$$

Where the elements of the information matrix $\hat{\mathbf{I}}_n(\hat{\boldsymbol{\theta}}) = \left(-\left(\frac{\partial^2 l(\boldsymbol{\theta})}{\partial\theta_i \partial\theta_j}\right)\right)_{\boldsymbol{\theta}=\hat{\boldsymbol{\theta}}}$ can be derived using the following second partial derivatives:

$$\frac{\partial^2 \ell}{\partial \theta^2} = -\frac{n}{\theta^2}$$

$$\frac{\partial^2 \ell}{\partial \alpha^2} = -\frac{n\theta}{\alpha^2} - (\theta+1)\sum_{i=0}^{n} \frac{[1-(1-e^{-\lambda t_i})^a]^{2b}}{\{1-\overline{\alpha}\,[1-(1-e^{-\lambda t_i})^a]^b\}^2}$$

$$\frac{\partial^2 \ell}{\partial a^2} = \frac{n}{a} + \sum_{i=0}^{n} \log(1-e^{-\lambda t_i})$$

$$+(1-b\theta)\sum_{i=0}^{n}\frac{(1-e^{-\lambda t_i})^{2a}\log(1-e^{-\lambda t_i})^2}{\{1-(1-e^{-\lambda t_i})^a\}^2}+(1-b\theta)\sum_{i=0}^{n}\frac{(1-e^{-\lambda t_i})^{a}\log(1-e^{-\lambda t_i})^2}{1-(1-e^{-\lambda t_i})^a}$$

$$+(\theta+1)\sum_{i=0}^{n}\frac{\bar{\alpha}^2 b^2[1-(1-e^{-\lambda t_i})^a]^{2(b-1)}(1-e^{-\lambda t_i})^{2a}\log(1-e^{-\lambda t_i})^2}{\{1-\bar{\alpha}[1-(1-e^{-\lambda t_i})^a]^b\}^2}$$

$$+(\theta+1)\sum_{i=0}^{n}\frac{b(b-1)\bar{\alpha}^2[1-(1-e^{-\lambda t_i})^a]^{b-2}(1-e^{-\lambda t_i})^{2a}\log(1-e^{-\lambda t_i})^2}{1-\bar{\alpha}[1-(1-e^{-\lambda t_i})^a]^b}$$

$$-(\theta+1)\sum_{i=0}^{n}\frac{b\,\bar{\alpha}\,[1-(1-e^{-\lambda t_i})^a]^{b-1}(1-e^{-\lambda t_i})^{a}\log(1-e^{-\lambda t_i})^2}{1-\bar{\alpha}[1-(1-e^{-\lambda t_i})^a]^b}$$

$$\frac{\partial^2\ell}{\partial b^2}=-\frac{n}{b^2}+(\theta+1)\sum_{i=0}^{n}\frac{\bar{\alpha}^2[1-(1-e^{-\lambda t_i})^a]^{2b}\log[1-(1-e^{-\lambda t_i})^a]^2}{\{1-\bar{\alpha}[1-(1-e^{-\lambda t_i})^a]^b\}^2}$$

$$+(\theta+1)\sum_{i=0}^{n}\frac{\bar{\alpha}\,[1-(1-e^{-\lambda t_i})^a]^{b}\log[1-(1-e^{-\lambda t_i})^a]^2}{1-\bar{\alpha}[1-(1-e^{-\lambda t_i})^a]^b}$$

$$\frac{\partial^2\ell}{\partial\lambda^2}=-\frac{n}{\lambda^2}+(a-1)\sum_{i=0}^{n}\left(-\frac{e^{-2\lambda t_i}t_i^{\,2}}{(1-e^{-\lambda t_i})^2}-\frac{e^{-\lambda t_i}t_i^{\,2}}{1-e^{-\lambda t_i}}\right)+(1-b\theta)\sum_{i=0}^{n}\frac{a^2(1-e^{-\lambda t_i})^{2(a-1)}e^{-2\lambda t_i}t_i^{\,2}}{\{1-(1-e^{-\lambda t_i})^a\}^2}$$

$$+(1-b\theta)\sum_{i=0}^{n}\frac{a(a-1)(1-e^{-\lambda t_i})^{a-2}e^{-2\lambda t_i}t_i^{\,2}}{1-(1-e^{-\lambda t_i})^a}-(1-b\theta)\sum_{i=0}^{n}\frac{a(1-e^{-\lambda t_i})^{a-1}e^{-\lambda t_i}t_i^{\,2}}{1-(1-e^{-\lambda t_i})^a}$$

$$+(\theta+1)\sum_{i=0}^{n}\frac{\bar{\alpha}^2a^2b^2e^{-2\lambda t_i}(1-e^{-\lambda t_i})^{2(a-1)}[1-(1-e^{-\lambda t_i})^a]^{2(b-1)}t_i^{\,2}}{\{1-\bar{\alpha}[1-(1-e^{-\lambda t_i})^a]^b\}^2}$$

$$+(\theta+1)\sum_{i=0}^{n}\frac{a^2\bar{\alpha}b(b-1)e^{-2\lambda t_i}(1-e^{-\lambda t_i})^{2(a-1)}[1-(1-e^{-\lambda t_i})^a]^{b-2}t_i^{\,2}}{1-\bar{\alpha}[1-(1-e^{-\lambda t_i})^a]^b}$$

$$-(\theta+1)\sum_{i=0}^{n}\frac{a(a-1)\bar{\alpha}\,b\,e^{-2\lambda t_i}(1-e^{-\lambda t_i})^{a-2}[1-(1-e^{-\lambda t_i})^a]^{b-1}t_i^{\,2}}{1-\bar{\alpha}[1-(1-e^{-\lambda t_i})^a]^b}$$

$$+(\theta+1)\sum_{i=0}^{n}\frac{a\,b\bar{\alpha}\,e^{-\lambda t_i}(1-e^{-\lambda t_i})^{a-1}[1-(1-e^{-\lambda t_i})^a]^{b-1}t_i^{\,2}}{1-\bar{\alpha}[1-(1-e^{-\lambda t_i})^a]^b}$$

$$\frac{\partial^2\ell}{\partial\theta\,\partial\alpha}=\frac{n}{\alpha}-\sum_{i=0}^{n}\frac{[1-(1-e^{-\lambda t_i})^a]^b}{1-\bar{\alpha}[1-(1-e^{-\lambda t_i})^a]^b}$$

$$\frac{\partial^2\ell}{\partial\theta\,\partial a}=-b\sum_{i=0}^{n}\frac{(1-e^{-\lambda t_i})^{a}\log(1-e^{-\lambda t_i})}{1-(1-e^{-\lambda t_i})^a}$$

$$-\sum_{i=0}^{n}\frac{b\bar{\alpha}[1-(1-e^{-\lambda t_i})^a]^{b-1}(1-e^{-\lambda t_i})^{a}\log(1-e^{-\lambda t_i})}{1-\bar{\alpha}[1-(1-e^{-\lambda t_i})^a]^b}$$

$$\frac{\partial^2\ell}{\partial\theta\,\partial b}=\sum_{i=0}^{n}\log[1-(1-e^{-\lambda t_i})^a]+\sum_{i=0}^{n}\frac{\bar{\alpha}[1-(1-e^{-\lambda t_i})^a]^b\log[1-(1-e^{-\lambda t_i})^a]}{1-\bar{\alpha}[1-(1-e^{-\lambda t_i})^a]^b}$$

$$\frac{\partial^2 \ell}{\partial \theta \partial \lambda} = -b\sum_{i=0}^{n} \frac{a\lambda (1-e^{-\lambda t_i})^{a-1} e^{-\lambda t_i}}{1-(1-e^{-\lambda t_i})^{a}} - \sum_{i=0}^{n} \frac{\overline{\alpha}\, a b \lambda\, [1-(1-e^{-\lambda t_i})^{a}]^{b-1} (1-e^{-\lambda t_i})^{a-1} e^{-\lambda t_i}}{1-\overline{\alpha}[1-(1-e^{-\lambda t_i})^{a}]^{b}}$$

$$\frac{\partial^2 \ell}{\partial \alpha \partial a} = (\theta+1) \sum_{i=0}^{n} \frac{\overline{\alpha}\, b\, [1-(1-e^{-\lambda t_i})^{a}]^{2b-1} (1-e^{-\lambda t_i})^{2a} \log(1-e^{-\lambda t_i})}{\{1-\overline{\alpha}[1-(1-e^{-\lambda t_i})^{a}]^{b}\}^{2}}$$

$$(\theta+1) \sum_{i=0}^{n} \frac{b\, [1-(1-e^{-\lambda t_i})^{a}]^{b-1} (1-e^{-\lambda t_i})^{a} \log(1-e^{-\lambda t_i})}{1-\overline{\alpha}[1-(1-e^{-\lambda t_i})^{a}]^{b}}$$

$$\frac{\partial^2 \ell}{\partial \alpha \partial b} = -(\theta+1) \sum_{i=0}^{n} \frac{\overline{\alpha}\, [1-(1-e^{-\lambda t_i})^{a}]^{2b} \log[1-(1-e^{-\lambda t_i})^{a}]}{\{1-\overline{\alpha}[1-(1-e^{-\lambda t_i})^{a}]^{b}\}^{2}}$$

$$-(\theta+1) \sum_{i=0}^{n} \frac{[1-(1-e^{-\lambda t_i})^{a}]^{b} \log[1-(1-e^{-\lambda t_i})^{a}]}{1-\overline{\alpha}[1-(1-e^{-\lambda t_i})^{a}]^{b}}$$

$$\frac{\partial^2 \ell}{\partial \alpha \partial \lambda} = (\theta+1) \sum_{i=0}^{n} \frac{\overline{\alpha}\, a b\, e^{-\lambda t_i} (1-e^{-\lambda t_i})^{a-1} [1-(1-e^{-\lambda t_i})^{a}]^{2b-1} t_i}{\{1-\overline{\alpha}[1-(1-e^{-\lambda t_i})^{a}]^{b}\}^{2}}$$

$$+(\theta+1) \sum_{i=0}^{n} \frac{a \overline{\alpha}\, b\, e^{-\lambda t_i} (1-e^{-\lambda t_i})^{a-1} [1-(1-e^{-\lambda t_i})^{a}]^{b-1} t_i}{1-\overline{\alpha}[1-(1-e^{-\lambda t_i})^{a}]^{b}}$$

$$\frac{\partial^2 \ell}{\partial a \partial b} = -\theta \sum_{i=0}^{n} \frac{(1-e^{-\lambda t_i})^{a} \log(1-e^{-\lambda t_i})}{1-(1-e^{-\lambda t_i})^{a}}$$

$$-(\theta+1) \sum_{i=0}^{n} \frac{\overline{\alpha} (1-e^{-\lambda t_i})^{a} [1-(1-e^{-\lambda t_i})^{a}]^{b-1} \log(1-e^{-\lambda t_i})}{1-\overline{\alpha}[1-(1-e^{-\lambda t_i})^{a}]^{b}}$$

$$-(\theta+1) \sum_{i=0}^{n} \frac{b \overline{\alpha}^{2} (1-e^{-\lambda t_i})^{a} [1-(1-e^{-\lambda t_i})^{a}]^{2b-1} \log(1-e^{-\lambda t_i}) \log[1-(1-e^{-\lambda t_i})^{a}]}{\{1-\overline{\alpha}[1-(1-e^{-\lambda t_i})^{a}]^{b}\}^{2}}$$

$$-(\theta+1) \sum_{i=0}^{n} \frac{b \overline{\alpha} (1-e^{-\lambda t_i})^{a} [1-(1-e^{-\lambda t_i})^{a}]^{b-1} \log(1-e^{-\lambda t_i}) \log[1-(1-e^{-\lambda t_i})^{a}]}{1-\overline{\alpha}[1-(1-e^{-\lambda t_i})^{a}]^{b}}$$

$$\frac{\partial^2 \ell}{\partial a \partial \lambda} = \sum_{i=0}^{n} \frac{e^{-\lambda t_i} t_i}{1-e^{-\lambda t_i}} + (1-b\theta) \sum_{i=0}^{n} \frac{a\, (1-e^{-\lambda t_i})^{2a-1} e^{-\lambda t_i} t_i \log(1-e^{-\lambda t_i})}{\{1-(1-e^{-\lambda t_i})^{a}\}^{2}}$$

$$+(1-b\theta) \sum_{i=0}^{n} \frac{(1-e^{-\lambda t_i})^{a-1} e^{-\lambda t_i} t_i}{1-(1-e^{-\lambda t_i})^{a}} + (1-b\theta) \sum_{i=0}^{n} \frac{a\, (1-e^{-\lambda t_i})^{a-1} e^{-\lambda t_i} t_i \log(1-e^{-\lambda t_i})}{1-(1-e^{-\lambda t_i})^{a}}$$

$$-(\theta+1) \sum_{i=0}^{n} \frac{b \overline{\alpha}\, e^{-\lambda t_i} (1-e^{-\lambda t_i})^{a-1} [1-(1-e^{-\lambda t_i})^{a}]^{b-1} t_i}{1-\overline{\alpha}[1-(1-e^{-\lambda t_i})^{a}]^{b}}$$

$$+(\theta+1) \sum_{i=0}^{n} \frac{a\, b^{2}\, \overline{\alpha}^{2} e^{-\lambda t_i} (1-e^{-\lambda t_i})^{2a-1} [1-(1-e^{-\lambda t_i})^{a}]^{2(b-1)} t_i \log(1-e^{-\lambda t_i})}{\{1-\overline{\alpha}[1-(1-e^{-\lambda t_i})^{a}]^{b}\}^{2}}$$

$$+(\theta+1) \sum_{i=0}^{n} \frac{a\, b (b-1) \overline{\alpha}\, e^{-\lambda t_i} (1-e^{-\lambda t_i})^{2a-1} [1-(1-e^{-\lambda t_i})^{a}]^{b-2} t_i \log(1-e^{-\lambda t_i})}{1-\overline{\alpha}[1-(1-e^{-\lambda t_i})^{a}]^{b}}$$

$$-(\theta+1)\sum_{i=0}^{n}\frac{a\,b\,\overline{\alpha}\,e^{-\lambda t_i}(1-e^{-\lambda t_i})^{a-1}[1-(1-e^{-\lambda t_i})^{a}]^{b-1}\,t_i\log(1-e^{-\lambda t_i})}{1-\overline{\alpha}[1-(1-e^{-\lambda t_i})^{a}]^{b}}$$

$$\frac{\partial^2\ell}{\partial b\partial\lambda}=-\sum_{i=0}^{n}\frac{a e^{-\lambda t_i}(1-e^{-\lambda t_i})^{a-1}t_i}{1-(1-e^{-\lambda t_i})^{a}}$$

$$-(\theta+1)\sum_{i=0}^{n}\frac{a\overline{\alpha}\;e^{-\lambda t_i}(1-e^{-\lambda t_i})^{a-1}[1-(1-e^{-\lambda t_i})^{a}]^{b-1}\,t_i}{1-\overline{\alpha}[1-(1-e^{-\lambda t_i})^{a}]^{b}}$$

$$-(\theta+1)\sum_{i=0}^{n}\frac{a\,b\overline{\alpha}^{2}\,e^{-\lambda t_i}(1-e^{-\lambda t_i})^{a-1}[1-(1-e^{-\lambda t_i})^{a}]^{2b-1}\,t_i\log[1-(1-e^{-\lambda t_i})^{a}]}{\{1-\overline{\alpha}[1-(1-e^{-\lambda t_i})^{a}]^{b}\}^{2}}$$

$$-(\theta+1)\sum_{i=0}^{n}\frac{a\,b\,\overline{\alpha}\,e^{-\lambda t_i}(1-e^{-\lambda t_i})^{a-1}[1-(1-e^{-\lambda t_i})^{a}]^{b-1}\,t_i\log[1-(1-e^{-\lambda t_i})^{a}]}{1-\overline{\alpha}[1-(1-e^{-\lambda t_i})^{a}]^{b}}$$

Where $\psi'(.)$ is the derivative of the digamma function.